\documentclass[a4paper,12pt,numbers,sort&compress]{article}

\pdfoutput=1

\usepackage[paper=letterpaper,margin=1in]{geometry}

\usepackage{amsmath,amssymb,amsthm,amsfonts,epsfig,cite,setspace,bigstrut,longtable,array,url,color}

\usepackage{hyperref}

\newcommand{\nn}{\nonumber}
\newcommand{\tr}{\mathop{\rm Tr}}
\newcommand{\comment}[1]{}

\newcommand{\cM}{{\cal M}}
\newcommand{\cW}{{\cal W}}
\newcommand{\cN}{{\cal N}}

\newcommand{\cO}{{\cal O}}

\newcommand{\cB}{{\cal B}}
\newcommand{\cC}{{\cal C}}
\newcommand{\cS}{{\cal S}}
\newcommand{\cD}{{\cal D}}
\newcommand{\cG}{{\cal G}}

\newcommand{\cQ}{{\cal Q}}

\newcommand{\IP}{\mathbb{P}}
\newcommand{\IQ}{\mathbb{Q}}

\newcommand{\IR}{\mathbb{R}}
\newcommand{\IC}{\mathbb{C}}
\newcommand{\IF}{\mathbb{F}}

\newcommand{\IZ}{\mathbb{Z}}

\newcommand{\drawsquare}[2]{\hbox{%
\rule{#2pt}{#1pt}\hskip-#2pt
\rule{#1pt}{#2pt}\hskip-#1pt
\rule[#1pt]{#1pt}{#2pt}}\rule[#1pt]{#2pt}{#2pt}\hskip-#2pt
\rule{#2pt}{#1pt}}
\newcommand{\fund}{\raisebox{-.5pt}{\drawsquare{6.5}{0.4}}}
\newcommand{\antifund}{\overline{\fund}}

\renewenvironment{thebibliography}[1]{%
\begin{oldthebibliography}{#1}%
\setlength{\parskip}{0ex}%
\setlength{\itemsep}{0ex}%
}%
{%
\end{oldthebibliography}%
}

\newcommand\blfootnote[1]{%
  \begingroup
  \renewcommand\thefootnote{}\footnote{#1}%
  \addtocounter{footnote}{-1}%
  \endgroup
}

\newtheorem{theorem}{\bf THEOREM}
\newtheorem{proposition}{\bf PROPOSITION}

\newtheorem{definition}{\bf DEFINITION}
\newtheorem{corollary}{\bf COROLLARY}
\newtheorem{observation}{\bf OBSERVATION}
\newtheorem{lemma}{\bf LEMMA}

\begin{document}
\begin{titlepage}

  ~\\
  \vspace{4mm}

\begin{center}
{\Large \bf  Calabi-Yau Varieties: from Quiver Representations to Dessins d'Enfants}
\medskip

\vspace{4mm}

{\large Yang-Hui He}\blfootnote{
  Invited chapter contribution to {\it Grothendieck-Teichm\"uller Theories and Impact},
  L.~Ji, A.~Papadopoulos, L.~Schneps, \& W.~Su Ed.~and based on talks given at the Chern Institute, Nankai University, Beijing Geometry Symposium, Perimeter Institute, Oberwolfach, and Harvard University.
}

\vspace{1mm}

\renewcommand{\arraystretch}{0.5} 
{\small
{\it
\begin{tabular}{rl}
  ${}^{1}$ &
  Department of Mathematics, City, University of London, EC1V 0HB, UK\\
  ${}^{2}$ &
  School of Physics, NanKai University, Tianjin, 300071, P.R.~China\\
  ${}^{3}$ &
  Merton College, University of Oxford, OX14JD, UK\\
\end{tabular}
}
~\\
~\\
~\\
hey@maths.ox.ac.uk
}
\renewcommand{\arraystretch}{1.5} 

\end{center}

\vspace{10mm}

\begin{abstract}
  The connections amongst (1) quivers whose representation varieties are Calabi-Yau, (2) the combinatorics of bipartite graphs on Riemann surfaces, and (3) the geometry of mirror symmetry have engendered a rich subject at whose heart is the physics of gauge/string theories.
  We review the various parts of this intricate story in some depth, for a mathematical audience without assumption of any knowledge of physics, emphasizing a plethora of results residing at the intersection between algebraic geometry, representation theory and number theory.
\end{abstract}

\end{titlepage}

\tableofcontents


\section{Introduction}
Of the myriad of colourful threads in the luxurious tapestry of mathematics there is a particular quartet which has caught the modern eye.
There is the representation theory of quivers, which has become a staple for the algebraist, and, since the McKay Correspondence, also that for an algebraic geometer, especially in the context of singularity resolution.
There is the geometry of Calabi-Yau varieties, which is itself a field of study, in being a paragon of K\"ahler and symplectic geometry, rich with such wonders as mirror symmetry.
There is the combinatorics of graphs endowed with bipartite structure, which is a corner stone to fields as varied as the statistical mechanics of dimer models to the geometry of positive Grassmannians.
And there is the arithmetic of the absolute Galois group, which, still elusive, is a central object to algebraic number theory.

There is, perhaps surprisingly, a strong skein which underlies these threads, and which weaves them seamlessly, and it comes from theoretical physics.
In particular, the study of supersymmetric gauge theories, notably in the context of string theory, has not only touched upon the aforementioned disciplines, but has inspired new results and directions of interest to both communities of mathematicians and physicists.
The subject has ranged from Calabi-Yau singularities and quiver algebras \cite{AU,beil,BKM,BCQ,rafCY,Broom,davison,DWZ,Ginz,ishii2010note,LM,MR,balazs}, to the physics of D-brane gauge theories in the AdS/CFT Correspondence \cite{Feng:2000mi,Feng:2001bn,Feng:2005gw}; from catalogues and computation of dessins \cite{catalog,hv,kmsv,sv,vk,Vidunas:2016xun} to their associated physical realization as brane-tilings \cite{Franco:2005rj,Franco:2005sm,Hanany:2005ve,gulotta,Davey:2009bp,Bose:2014lea,Hanany:2011ra,Hanany:2011bs,He:2012xw,He:2014jva,He:2012jn,He:2015vua,Jejjala:2010vb}; from integrable models \cite{GK,Eager:2011dp} and mirror symmetry \cite{Benvenuti:2006qr,Franco:2012wv,Franco:2016qxh,He:2016yvz,Heckman:2012jh,Hanany:2008sb,Martelli:2006yb} to the representation theory of cluster algebras \cite{FZ,Mar,BKM}, as well as to curious infringements upon number theory \cite{He:2016yvz,Zhou:2015tia,He:2012jn,He:2012xw}.

Whilst there are excellent reviews \cite{Kennaway:2007tq,Yamazaki:2008bt} (cf.~also a rapid synopsis in \cite{He:2012js}), that these were aimed at string theorists and that the field has since advanced rapidly, beckon for a somewhat new comprehensive overview for a pure mathematician.
The thesis \cite{Broom} is a marvelous source but again the length and its aim toward algebraic geometers might deter a more general audience.
One is thus confronted with a formidable challenge: to have a medium-portioned taster to a mathematician, enough to more than merely whet an appetite but not so much as to overwhelm a glutton.

To meet this challenge we will certainly fail.
Yet, given this opportunity to speak to a diverse mathematical audience at the Chern Institute in this excellent conference organized by Professors L.~Ji, A.~Papadopoulos, L.~Schneps, \& W.~Su on Grothendieck-Teichm\"uller Theories, we feel compelled to at least attempt something appropriate for a post-graduate student in mathematics who wishes for a speedy entry into so vast and wondrous a domain and we hope that this chapter will be of some use.

Though we will be maintaining a dialogue between the mathematics and the physics by keeping a dictionary of terminology, the definitions, examples and propositions throughout will, in as much clarity as we can manage, be structured for a reader with no knowledge of physics whatsoever.
We begin with the definition of a quiver $\cQ$ and its associated representation variety $\cM$ in \S\ref{s:quiver}, emphasizing when $\cM$ is affine Calabi-Yau.
Next, in \S\ref{s:bi}, we show how in the case of $\cM$ being toric, $\cQ$ is dual to a bipartite graph drawn on the torus, and how cluster transformation on the quiver exhibits in various facets of the story.
We then discuss in \S\ref{s:dessin} how the shape of the torus is determined from the point of view of geometry and from that of dessins d'enfants and what puzzles arise.
Finally, we conclude with open questions in \S\ref{s:conc}.

\section{Quivers and Moduli Spaces}\label{s:quiver}
Our starting point is a quiver with superpotential, with the following data:
\begin{definition}
  A quiver $\cQ = (\cQ_0, \cQ_1, W)$ is a finite directed graph with the set of vertices $\cQ_0$ and arrows $\cQ_1$, the cardinalities of which are $N_0$ and $N_1$ respectively:
  \begin{itemize}
  \item $\cQ$ is equipped with a {\em representation}, meaning that we attach $V_i \simeq \IC^{n_i}$ to each node for some positive integer $n_i$, whence each arrow $(X_{ij} \in \cQ_1) \in \hom(V_i, V_j)$ can be considered as an $n_j \times n_i$ matrix; we allow self-adjoining arrows, $\phi_i = X_{ii}$, as well as cycles which are closed loops $X_{i_1i_2} X_{i_2i_3} \ldots X_{i_ki_1}$ formed by the arrows.
  \item $\cQ$ is also furnished by {\em relations}; this is imposed by the superpotential $W$ which is a polynomial in all the arrows treated as formal matrix variables:
    \begin{equation}\label{W}
    W = \sum_{k=1}^{N_2} c_k \tr(\prod X_{ij}) \ldots \tr(\prod X_{i'j'})
    \end{equation}
    summed over possible cycles or products therein with coefficients $c_k \in \IC$. The formal polynomial relations amongst the arrows are determined by the vanishing of the Jacobian $\partial_{X_{ij}} W$.
  \end{itemize}
\end{definition}

\paragraph{Remark }
We supplant the above definition with several remarks
\begin{itemize}
\item The allowance for loops and cycles significantly complicates the representation theory of $\cQ$ but, as we shall shortly see, this is necessary for the physics.
  
\item The numbers $N_0$, $N_1$ and $N_2$ (the total number of monomial terms in $W$) will play an interesting combinatorial r\^ole.
  
\item Though a centre-piece of representation theory and algebraic geometry, $\cQ$ has also, perhaps unexpectedly, become a indispensable part of modern theoretical physics.
  Into the details of quantum field theories we shall not delve and it suffices to say that to the above quiver data is associated a 4-dimensional, supersymmetric $\cG = \prod\limits_{i=1}^{N_0} U(n_i)$ gauge theory, under the following dictionary:
  
  {\small \hspace{-1cm}
  \begin{tabular}{ccc}
    Node $i$ && factor $U(n_i)$ \\
    Arrow $i \to j$ && bi-fundamental field $X_{ij}$ transforming as
    $(\fund,\antifund)$ of $U(N_i)\times U(N_j)$ \\
    Loop $i \to i$ && adjoint field $\phi_i = X_{ii}$ of $U(N_i)$\\
    Cycle $i_1 \to i_1 \to \ldots i_k \to i_1$ && Gauge Invariant Operator (GIO)
    $\tr(X_{i_1i_2} X_{i_2i_3} \ldots X_{i_ki_1})$ \\
    $\tr(\prod X_{ij})$ && Single-trace GIO\\
    $\tr(\prod X_{ij}) \ldots \tr(\prod X_{i'j'})$ && Multi-trace GIO\\
    2-Cycle $X_{ij}X_{jk}$ && Mass term\\
    Superpotential \eqref{W} && Superpotential in the Lagrangian with couplings $c_i$\\
    $\left\{ \partial_{X_{ij}} W \right\}$ && F-Terms\\
  \end{tabular}
  }
  
  The mathematician so inspired to converse with a physicist can use this dictionary though we shall mostly adhere to the left column for this review.
  
\item The list of labels $\vec{n} = (n_1, n_2, \ldots, n_{N_0})$ is called the dimension vector and for this note, we will take, for convenience, all $n_i = 1$ so that the arrows are just complex numbers; our quantum field theory thus has gauge group $U(1)^{N_0}$.

\item Suppose the {\em incidence matrix} of $\cQ$ is $d_{i\alpha}$ where $i = 1, \ldots , N_0$ indexes the nodes and $\alpha = 1, \ldots, N_1$, the arrows, such that each arrow $i \to j$ gives a column in $d_{i\alpha}$ with $-1$ at row $i$ and $+1$ at row $j$, and 0 otherwise.
  Then $\sum\limits_{\alpha} d_{j \alpha} |X_{ij}|^2 - \zeta_i$ are known as {\it D-terms}, where $\zeta_i \in \IC$ are so-called Fayet-Iliopoulos (FI) parametres \footnote{
    In general, one assigns charges $q_\alpha$ to the fields and sum over $q_\alpha |X_{i\alpha}|^2$ but for our present purposes of $U(1)^{N_0}$ theories, the incidence matrix serves to encode the charges.
    Moreover, the FI-parametres exist only for the $U(1)$ gauge group factors.
  }.
  Note that unlike the F-terms, these are non-holomorphic, as they involve the complex conjugates of $X_{ij}$.
\end{itemize}

One of the most important quantities associated to a quantum field theory is its vacuum; in supersymmetric theories, the vacuum is often a continuous manifold parametrized by the expectation values of the scalar field components of the various multiplets.
From the quiver perspective, this can be realized rather succinctly \cite{naka,king}.
\begin{definition}
  The quiver variety $\cM(\cQ)$  is the GIT quotient of the representations $Rep(\cQ) = \oplus_{i,j} \hom(\IC^{n_i}, \IC^{n_j})$, with relations from $W$, quotiented by the complexified group $\cG_{\IC} = \prod_i GL_n(\IC)$
  \[
  \cM(\cQ) = {\rm Spec} \IC[ Rep(\cQ) / \left< \partial_{X_{ij}} W \right>]^{\cG_{\cC}} \simeq \left< \partial_{X_{ij}} W \right> // \cG_{\IC} \ .
  \]
  In physics, $\cM(\cQ)$ is called the {\em vacuum moduli space} (VMS) of the gauge theory \footnote{
    Strictly speaking, $\cM$ here defined is the mesonic branch of the moduli space because the gauge invariants are built from taking traces. We could contract with other invariant tensors, for instance, contracting with Levi-Civita symbols to obtain the GIOs would give the baryonic moduli space.
    }.
\end{definition}
From the point of view of geometric representation, $\cM(\cQ)$ can be construed as the centre of the path algebra of $\cQ$.
The Jacobian variety $\left< \partial_{X_{ij}} W \right>$ is itself also interesting and is called the {\em master space} \cite{ot,Forcella:2008bb} and the equations $\partial_{X_{ij}} W = 0$ are collectively called the {\it F-terms}.

The VMS is a supersymmetric vacuum manifold of the complex scalars, and is explicitly the solution set to the vanishing of the F-terms and the D-terms, known as F-flatness and D-flatness:  $\cM = \{F=D=0\}$.
The non-holomorphic nature of the D-terms render the computation of $\cM$ rather difficult and one resorts to the above GIT quotient, which though seemingly abstract, actually gives a very explicit way \cite{Mehta:2012wk,book} of obtaining $\cM(\cQ)$:
\begin{proposition}\label{vms}
  The VMS $\cM(\cQ)$ is realized as an affine algebraic variety by the following algorithm
  \begin{enumerate}
  \item Let $GIO_{min}$ be the set of minimal (Eulerian) cycles in $\cQ$, i.e., $\Phi_{r=1,\ldots,k}$, each of which is a polynomial in $X_{ij}$ and $k$, the total number of such cycles;
  \item Consider the ring map \footnote{
    As coordinate rings the map should go in the reverse direction so that as varieties the map goes as indicated, but we beg the readers' momentary indulgence so that it is clear that the {\it image} is an affine variety in $\IC^k$.
    } from the quotient ring by the Jacobian ideal by $\Phi_r$ as 
    $$\IC[X_{ij}] / \left< \partial W \right> \longrightarrow \IC[\phi_r] \ ; $$
  \item $\cM$ is the image of this map, as an affine variety in $\IC^k$.
  \end{enumerate}
\end{proposition}

\subsection{Example}
It is illustrative to pause a while for a concrete example.
Consider the $U(1)^3$ quiver with superpotential $W$ (where $\epsilon$ is the totally antisymmetric symbol on 3 indices) as follows
\begin{equation}\label{egC3Z3}
\begin{array}{cc}
\begin{array}{c} \includegraphics[trim= 0mm 0mm 0mm 150mm,clip,width=2.5in]{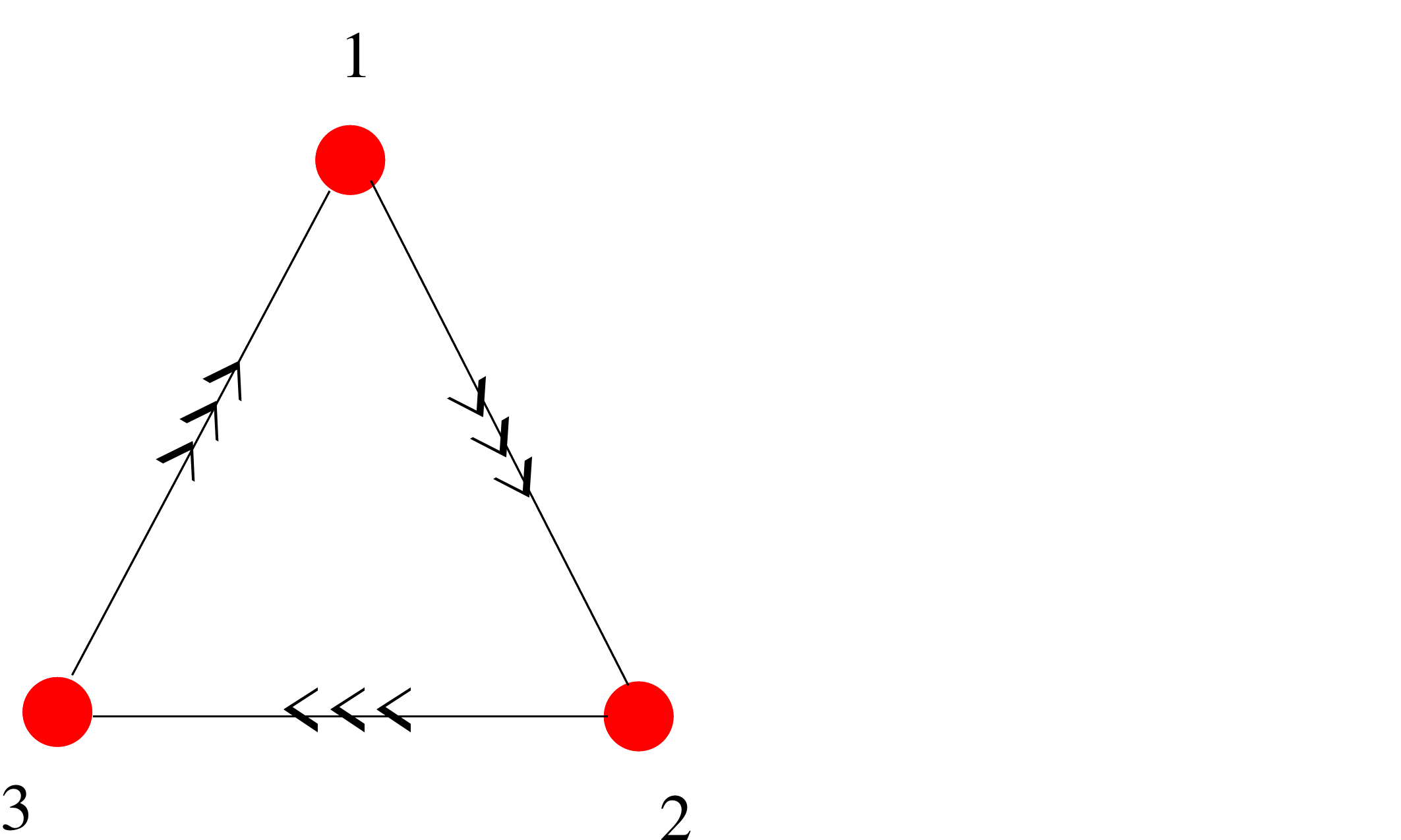} \end{array}
&
W= \sum\limits_{\alpha,\beta,\gamma=1}^3\epsilon_{\alpha\beta\gamma} X^{(\alpha)}_{12} X^{(\beta)}_{23} X^{(\gamma)}_{31}.
\end{array}
\end{equation}
There are nine arrows (fields) $X^{(\alpha)}_{12}, X^{(\beta)}_{23}$, and $X^{(\gamma)}_{31}$, $\alpha,\beta,\gamma=1,2,3$, with the subscript $ij$ signifying an arrow from node $i$ to $j$ and superscript $\alpha$ denoting that there is a multiplicity of three arrows for each pair of nodes.
Hence, $N_0 = 3$, $N_1 = 9$, and expanding out the $\epsilon$, $N_2 = 6$.
The Eulerian cycles are clearly the $3^3 = 27$ GIOs formed by their products corresponding to the closed cycles in the quiver.
Note that because the dimension vector is chosen to be $(1,1,1)$ all 9 arrows are just complex numbers (as elements of $\hom(\IC,\IC)$).

Taking the partial derivatives of $W$ with respect to the fields, we obtain 9 F-terms: $ \sum\limits_{\beta,\gamma=1}^3\epsilon_{\alpha\beta\gamma} X^{(\beta)}_{23}X^{(\gamma)}_{31}$,
$\sum\limits_{\alpha,\gamma=1}^3\epsilon_{\alpha\beta\gamma} X^{(\alpha)}_{12}X^{(\gamma)}_{31}$, and
$\sum\limits_{\alpha,\beta=1}^3\epsilon_{\alpha\beta\gamma} X^{(\alpha)}_{12}X^{(\beta)}_{23}$.
The VMS is thus
\begin{equation}
\cM \simeq {\rm Im} \left( \frac{\IC[X^{(\alpha)}_{12}, X^{(\beta)}_{23},X^{(\gamma)}_{31}]}
                                {\langle \epsilon_{\alpha\beta\gamma} X^{(\beta)}_{23}X^{(\gamma)}_{31},
                                                \epsilon_{\alpha\beta\gamma} X^{(\alpha)}_{12}X^{(\gamma)}_{31},
                                                \epsilon_{\alpha\beta\gamma} X^{(\alpha)}_{12}X^{(\beta)}_{23} \rangle}
                           \stackrel{GIO = \{X^{(\alpha)}_{12} X^{(\beta)}_{23}X^{(\gamma)}_{31} \}}{\xrightarrow{\hspace*{3cm}}}
                           \IC[\Phi_1, \ldots, \Phi_{27}] \right) \ ,
\end{equation}
where in the Jacobian ideal of $W$, we for simplicity of notation, use the summation convention that repeated indices are summed.
We find, using for example \cite{mac}, that $\cM$ is an affine variety of complex dimension three and degree nine;
it is the intersection of $17$ lines and $27$ quadratics in $\IC^{27}$.
The astute reader might recognize these as the defining polynomials for the total space of the anti-canonical line bundle $\cO_{\IP^2}(-3)$ and to this point we now turn.

\subsection{Calabi-Yau VMS}
In general, given a quiver $\cQ$, $\cM(\cQ)$ could be any affine variety (physically, this means that the VMS of a supersymmetric QFT could potentially be of any shape and form).
However, the vast majority of the literature is concerned with when $\cM(\cQ)$ is an affine Calabi-Yau variety of dimension 3.
The motivation of this originates from string theory.

Briefly, there are objects known as D3-branes in the theory which are $3+1$-dimensional subspaces of the $\IR^{9,1}$ background of the superstring.
The world-volume of a dynamical D3-brane supplies a 4-dimensional supersymmetric gauge theory whose VMS, by construction, is the transverse $10-3=6$ dimensions.
In order to preserve supersymmetry, these transverse 6 directions can be chosen to be a Ricci-flat complex (K\"ahler) 3-fold, i.e., an affine Calabi-Yau threefold (CY3).
Thus a natural correspondence is established between the physics of the gauge theory and the complex geometry of the affine CY3; this can be conceived of as an algebro-geometric version of the famous AdS/CFT (holographic) Correspondence.

Henceforth, we will focus on $\cM(\cQ)$ being an affine CY3, which, for our purposes, can be taken as a variety with trivial canonical sheaf (and admitting a crepant resolution to a smooth CY3).
For instance, the above case of $\cO_{\IP^2}(-3)$ is a perfect example and locally, this is a $\IZ/3\IZ$ orbifold of $\IC^3$.

The most archetypal example, however, is indubitably the ``clover quiver''
\begin{equation}\label{c3}
  \begin{array}{c}
    \includegraphics[trim=0mm 0mm 0mm 0mm, clip, width=1in]{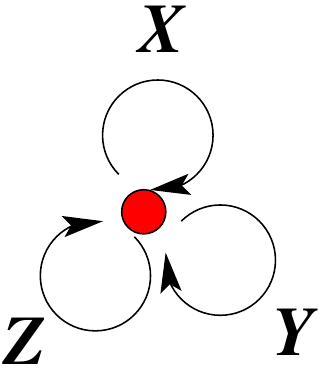}
  \end{array}
  \qquad
  W = \tr(XYZ - XZY) \ ,
\end{equation}
where we momentarily restore the label to the node to be $n$.
Clearly, the F-terms demand that $[X,Y]=[Y,Z]=[X,Z]$ and whence we have $\cM$ parametrized by 3 mutually commuting $n \times n$ matrices, which is nothing other than the $n$-th symmetric product of $\IC^3$.
That is to say $\cM \simeq (\IC^3)^n / S_n$.
For $n=1$, the VMS is simply $\IC^3$, the simplest affine CY3.
The gauge theory encoded by this data is known as $\cN=4$ super-Yang-Mills theory in 4-dimensions, perhaps the most well-known of all supersymmetric quantum field theories.

Whilst no Calabi-Yau metrics are known for compact smooth cases, those for affine CY3 are now known in abundance;
these are all realized as conical metrics over Sasaki-Einstein 5-manifolds.
Furthermore, over the last two decades there has been a host of activity in both the physics and the mathematics communities to construct explicit affine (or local) CY3s.
These have included quotienting $\IC^3$ by discrete finite subgroups of $SU(3)$, anti-canonical bundles over del Pezzo surfaces and affine toric varieties.
We summarize the foregoing discussions in the following diagram:
\[
\begin{array}{cc}
  \begin{array}{l}\includegraphics[width=2.5in,angle=0]{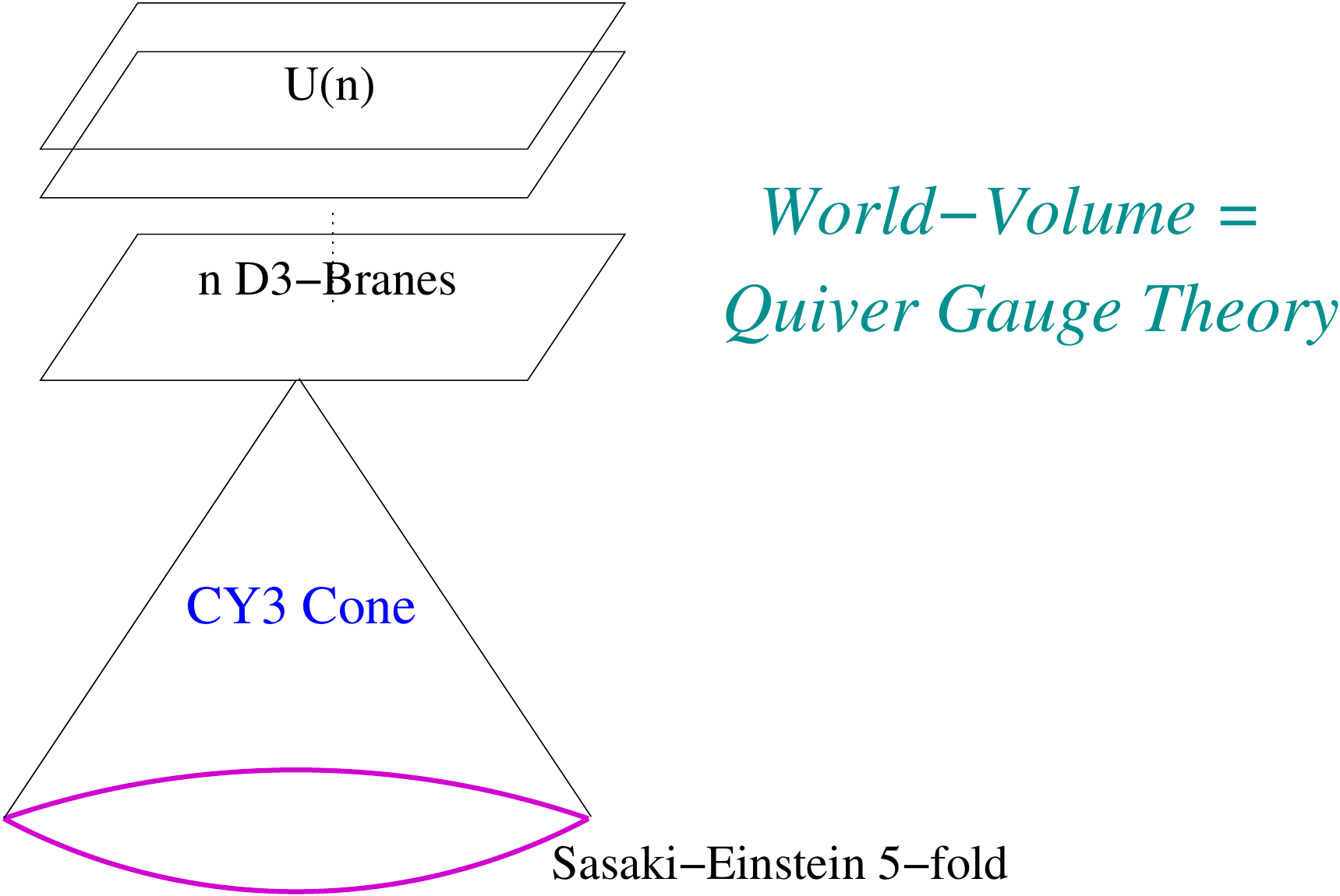}\end{array}
  &
  \begin{array}{l}\includegraphics[width=3in,angle=0]{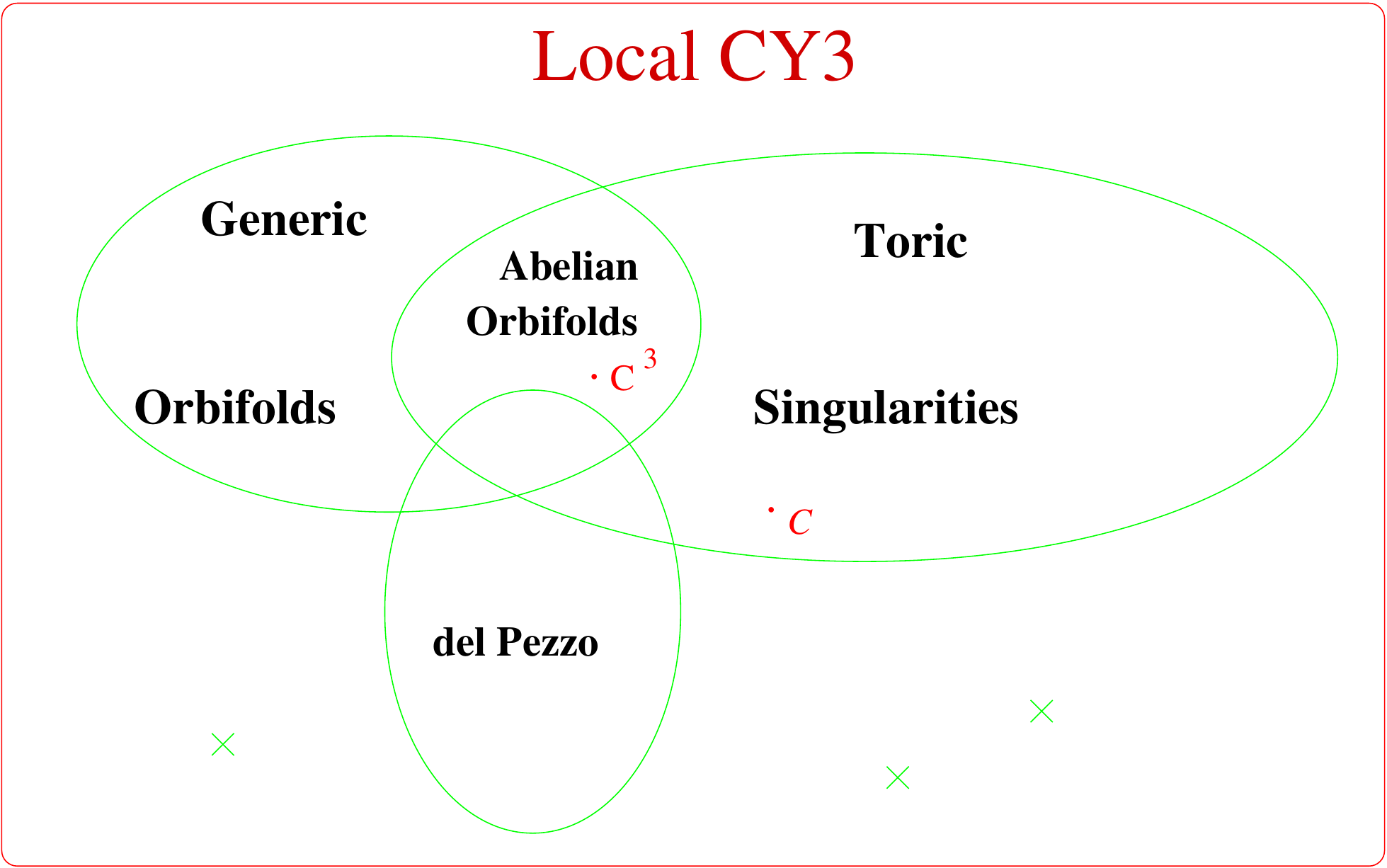}\end{array}
\end{array}
\]
On the left, we have sketched the physical situation of D3-branes in string theory placed at the tip of a Calabi-Yau cone over a Sasaki-Einstein five-manifold; the world-volume of the branes is a gauge theory encoded by the quiver $\cQ$, whose VMS is, by construction, $\cM(\cQ)$.
On the right we give the state of the art knowledge of the space of affine Calabi-Yau threefolds.
In the Venn diagram, we have marked the two most famous local CY3s: $\IC^3$ and $\cC$, the quadric $uv=wz$ in $\IC[u,v,w,z]$ (known as the conifold in physics).
The green crosses outside the main ellipses are isolated constructions.
It should be noted that the quiver for the orbifold case is none other than the generalized McKay quiver \cite{Hanany:1998sd} and also the toric singularities are by far (consisting of many infinite families) the largest and the best studied cases; to these we now turn.

\subsubsection{Affine Toric Calabi-Yau 3-folds}\label{s:CY3}
We briefly remind the reader the key facts for affine toric varieties.
Suppose we have a (strongly convex rational polyhedral) cone $\sigma \subset \IZ^r \otimes \IR$, then $X = $Spec$\IC[\sigma^\vee \cap \IZ^r]$ is an affine toric variety of dimension $r$ associated with the cone $\sigma$; here $\sigma^\vee$ is the dual cone
\footnote{Strictly speaking we should be taking a pair of dual integer lattices $N,M \simeq \IZ^r$ so that $\sigma \subset N \otimes_{\IZ} \IR$ and $\sigma \subset M \otimes_{\IZ} \IR$, and we computing $\sigma^\vee \cap M$, but we will be cavalier and take advantage of the isomorphism between $\IZ^r$ and its dual to just use $\IZ^r$.}.
The toric diagram $\cD$ is simply the end points of the generators of $\sigma$; this is simply a list of integer $r$-vectors.

Importantly, the Calabi-Yau condition of having trivial anti-canonical sheaf demands that all the vectors are co-hyperplanar, so that we can draw $\cD$ in as a grid of lattice points in $\IZ^{r-1}$.
Now, we are interested in 3-folds, so $r=3$, and we have the rather convenient
\begin{proposition}
Every convex lattice polygon gives an affine toric Calabi-Yau threefold.
\end{proposition}

For instance, the cone for $\IC^3$ is generated by the standard vectors $(1,0,0)$, $(0,1,0)$ and $(0,0,1)$, so the toric diagram is the triangle as given in part (a) of the following figure.
Likewise, the quadric $\cC$ is also toric Calabi-Yau whose $\cD$ is the square as given in part (b).
\begin{equation}\label{toricdiag}
  (a) \ 
  \IC^3
  \begin{array}{c}
    \includegraphics[trim= 0mm 0mm 0mm 0mm,clip,width=1.5in]{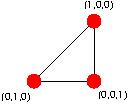}
  \end{array}
  \qquad \qquad
  (b) \
  \cC
  \begin{array}{c}
    \includegraphics[trim= 0mm 0mm 0mm 0mm,clip,width=1.5in]{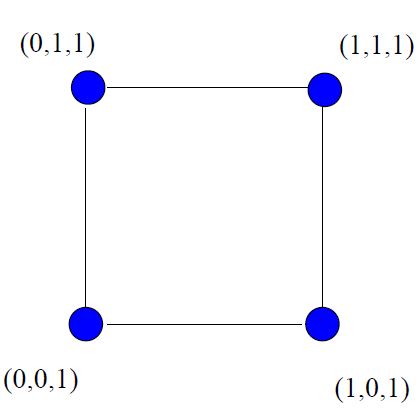}
  \end{array}
\end{equation}
We see, therefore, why there is a plethora of affine toric CY3 \footnote{Note that no compact CY3 is toric.}; for many of these explicit metrics have been found (cf,~\cite{Martelli:2006yb}).
Several typical examples have enriched the literature and it is expedient for us to summarize these into a catalogue, especially since some of the appellations are rather esoteric (cf.~\cite{Closset:2009sv}).
In the ensuing we use physicists' notation $\IZ_k$ to mean the cyclic group $\IZ / k \IZ$.

\newpage

\begin{longtable}{>{\centering}m{6cm} |c}
  Toric Diagram & Description \\
  \hline
  $\begin{array}{c}
    \includegraphics[trim= 0mm 0mm 0mm 0mm,clip,width=1.5in]{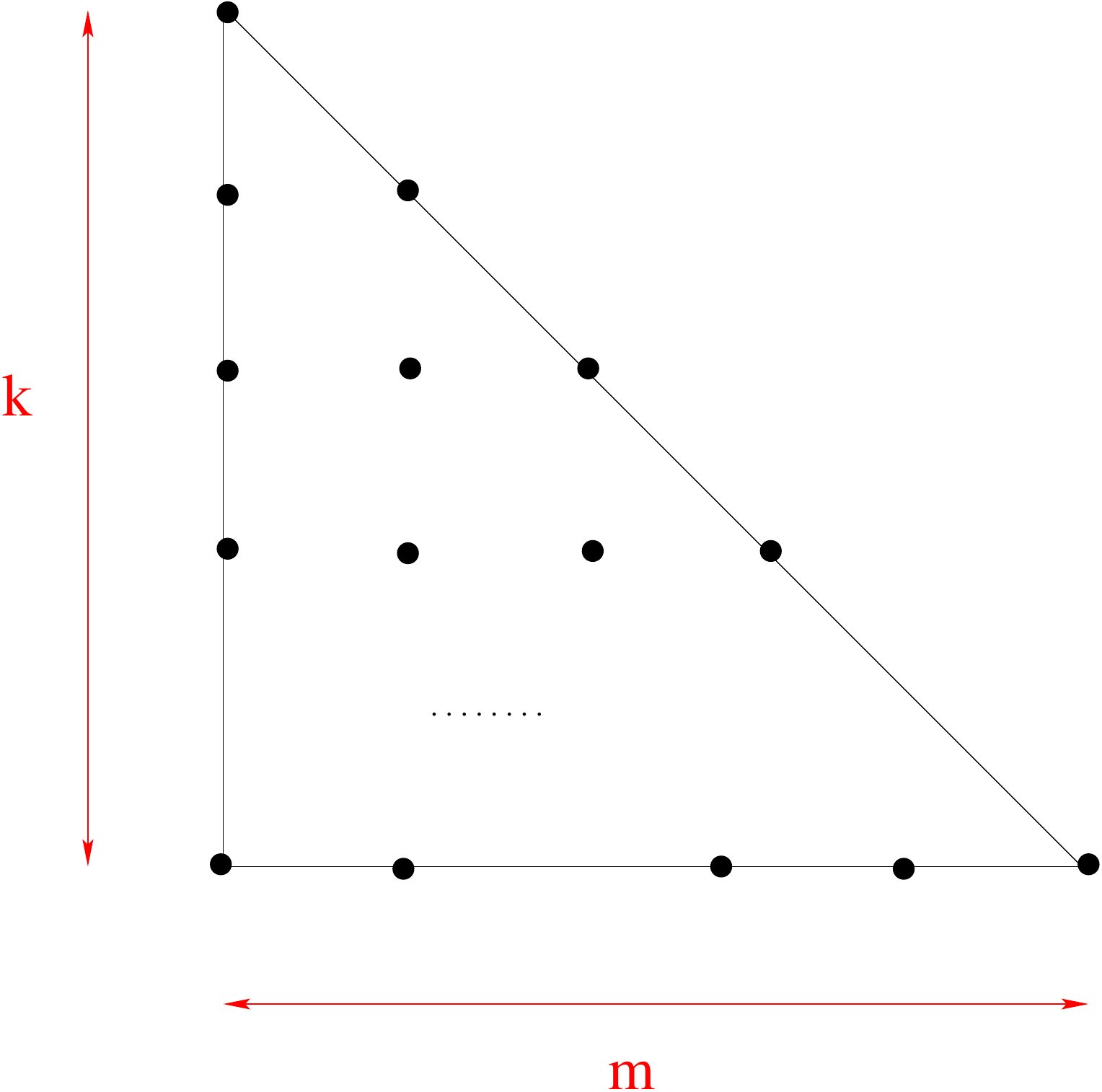}
  \end{array}$
  &
  \begin{tabular}{l}
    {\bf Abelian orbifold $\IC^3 / (\IZ_k \times \IZ_m)$}\\
    with action by the generators of the 2 cyclic factors as:\\
    $\IC^3 \ni (x,y,z) \mapsto (\omega_k x, \omega_k y, \omega_k^{-2} z) ; \
    (\omega_m x, \omega_m y, \omega_m^{-2} z)$ \\
    $(k,m)=(1,1)$ is simply $\IC^3$ \\
  \end{tabular} 
  \\    \hline
  $\begin{array}{c}
    \includegraphics[trim= 0mm 0mm 0mm 0mm,clip,width=2in]{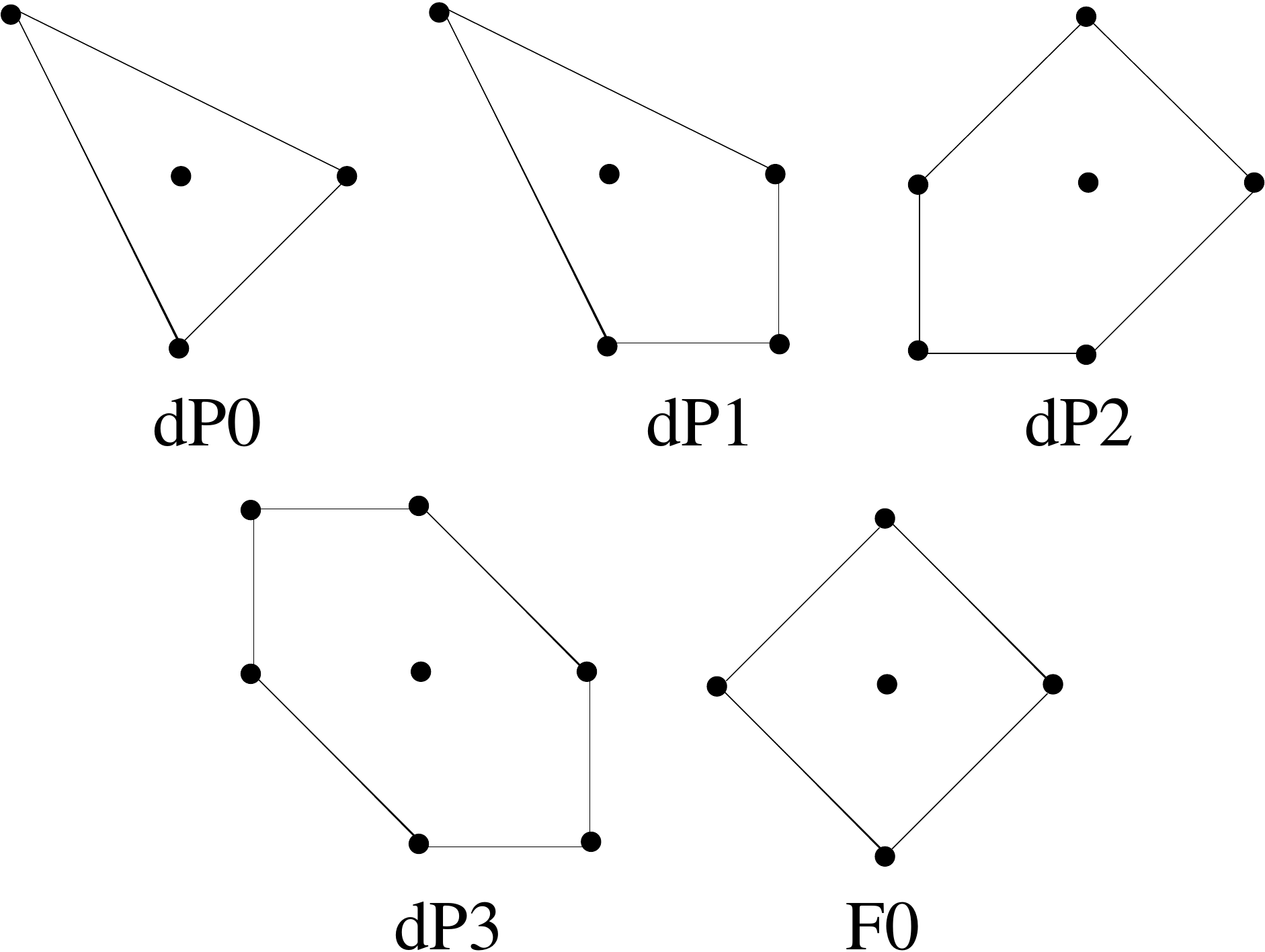}
  \end{array}$
  &
  \begin{tabular}{l}
    {\bf complex cones over toric Fano surfaces}\\
    $dP_n$ is del Pezzo surface: $\IP^2$ blown up at $n$ generic points,\\
    \quad $n$ can go up to 8 but only first 4 are toric; \\
    $dP_0 \simeq \IC^3 / \IZ_3$, \ 
    $F_0 \simeq \IP^1 \times \IP^1$ is zeroth Hirzebruch surface
  \end{tabular} 
  \\    \hline
  $\begin{array}{c}
    \includegraphics[trim= 0mm 0mm 0mm 0mm,clip,width=2in]{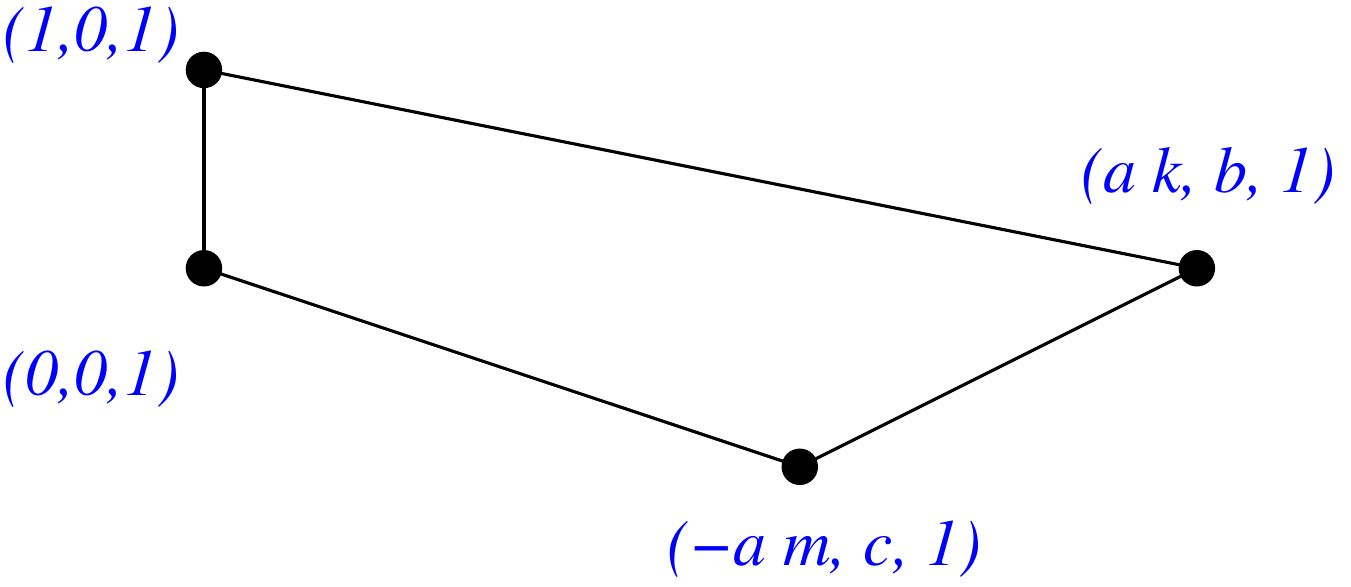}
  \end{array}$
  &
  \begin{tabular}{l}
    {\bf $L^{a,b,c}$}, \quad $a,c \leq b \in \IZ_+$, \ $\gcd(a,b,c,a+b-c) = 1$ \\
    In the figure, $c k + b m = 1$ \\
    Special cases: $L^{p-q,p+q,p} = Y^{p,q}$, $L^{1,1,1} = \cC$ (conifold),\\
    \qquad $L^{1,2,1} = SPP$ (suspended pinch point)
  \end{tabular}
  \\    \hline
  $\begin{array}{c}
    \includegraphics[trim= 0mm 0mm 0mm 0mm,clip,width=2in]{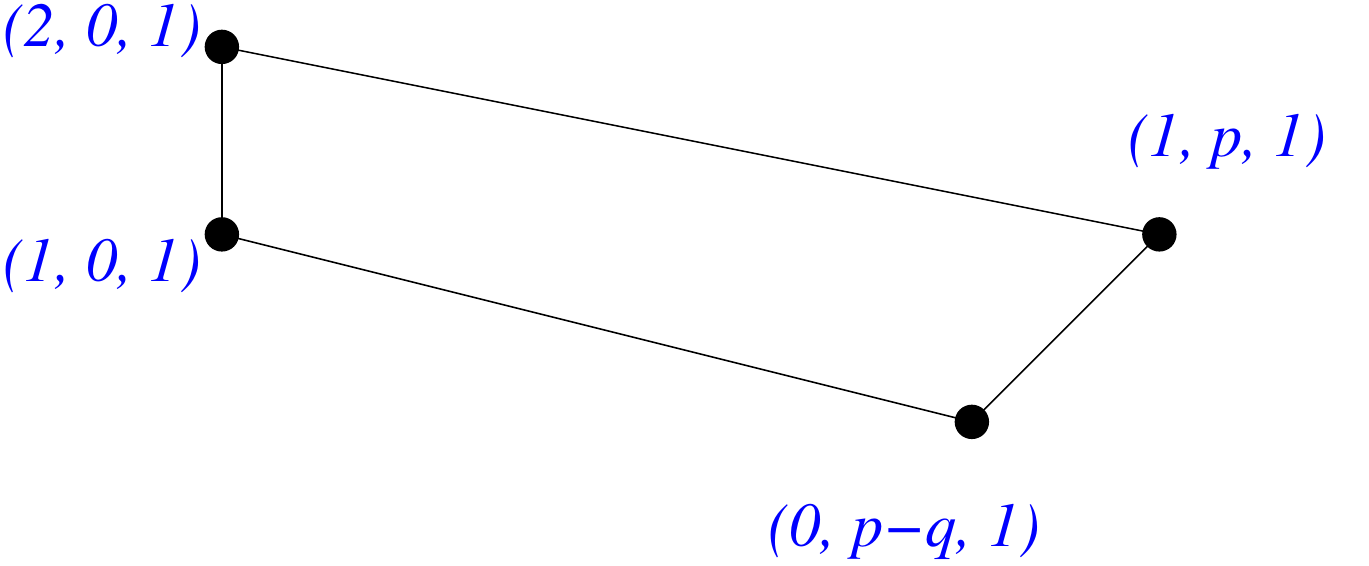}
  \end{array}$
  &
  \begin{tabular}{l}
    {\bf $Y^{p,q}$}, \quad $q < p \in \IZ_{\geq 0}$\\
    Topologically, real cone over $S^2 \times S^3$; \\
    $Y^{2,1} = dP_1$; \quad
    $Y^{p,0} = \cC / \IZ_p$ \\
    Extended definition: $Y^{p,p} = \IC^3 / \IZ_{2p}$\\
  \end{tabular}
  \\    \hline
  \\[0.05cm]
  $\begin{array}{c}
    \includegraphics[width=2in]{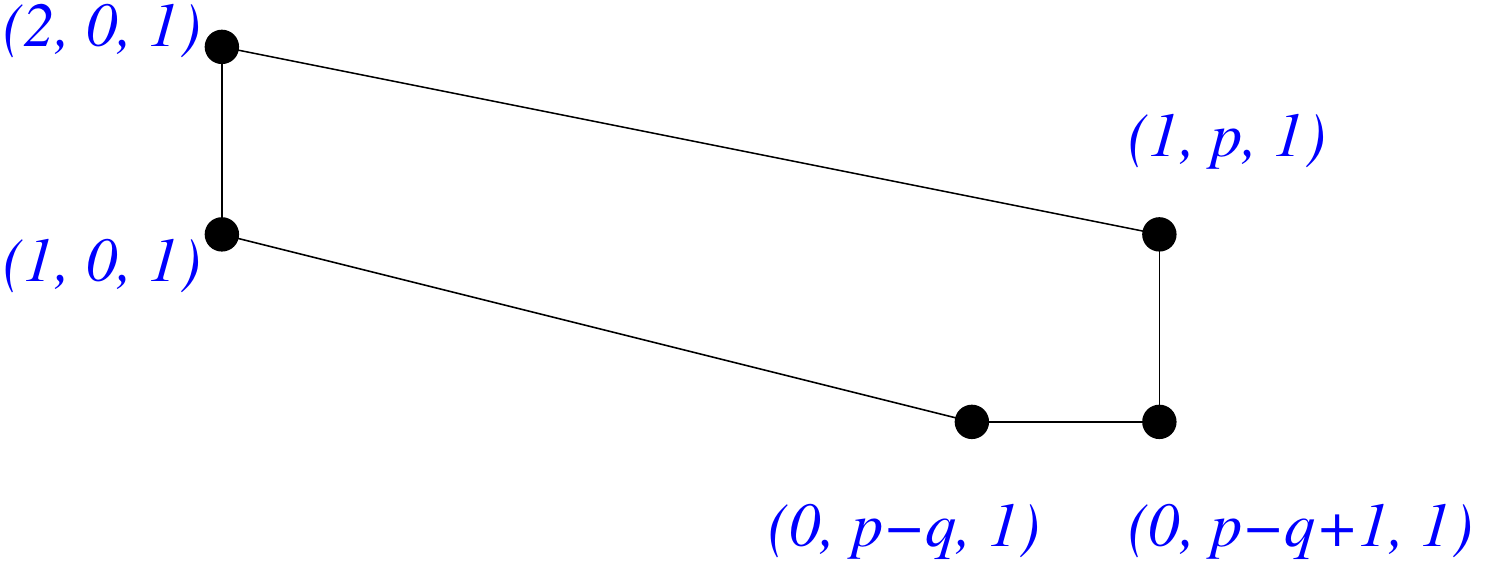}
  \end{array}$
  &
  \begin{tabular}{l}
    {\bf $X^{p,q}$},\quad $q < p \in \IZ_{\geq 0}$\\
    Special cases: $X^{1,1} = SPP$, \ $X^{2,1} = dP_2$\\
  \end{tabular}
\end{longtable}

As complex varieties, the above are affine toric CY3 $\cM$, and it is an important fact (especially to the physics) that they are all real cones over some Sasaki-Einstein 5-manifolds (SE5) whose metrics have been explicitly constructed in many cases, the archetypal one is perhaps $\cC = L^{1,1,1}$ whose SE5 is the so-called $T^{1,1}$ space.

Because all our toric diagrams are planar, they are easily encoded into a single bivariate polynomial:
\begin{definition}\label{newton}
  Let the coordinates of all (including boundary and interior) lattice points of the toric diagram $\cD$ be $(a_i,b_i)$ and $(z,w) \in \IC^2$ be formal variables, the Laurent polynomial
  \[
  P_{\cD}(z,w) = \sum_{(u_i, v_i) \in \cD} \alpha_i z^{a_i} w^{b_i}
  \]
  is called the Newton polynomial of $\cD$.
\end{definition}
In the above, $\alpha_i$ are complex coefficients which are parametres.
Although $P_{\cD}$ is a Laurent polynomial since the coordinates of the toric diagram can be negative, since $\cD$ is defined only up to shifts (and $SL(2;\IZ)$) in $\IZ^2$, we can always make it into a true polynomial.
It is a central theorem of local mirror symmetry \cite{Strominger:1996it,Hori:2000kt} that
\begin{theorem}
  The affine variety $\cW$ defined by
  \[
  \{ uv = P_{\cD}(z,w) \} \subset \IC^4[u,v,z,w]
  \]
  is the local mirror to $\cM$ with toric diagram $\cD$.
\end{theorem}
We note that though $\cW$ is affine CY3, it is not necessarily toric.
The parametres $\alpha_i$ manifest as complex structure of $\cW$.

\section{Bipartite Graphs on $T^2$}\label{s:bi}
Returning to our examples presented in \eqref{egC3Z3} and \eqref{c3} at the beginning of the previous section, we notice two things.
First, there is a curious combinatorial relation
\begin{equation}\label{T2}
N_0 - N_1 + N_2 = 0
\end{equation}
and second, that all our VMS are toric.
This turns out not to be a coincidence.

The origin of the condition \eqref{T2} is the superconformal symmetry of our gauge theories in the infra-red, though we only mention this in passing (q.v.~\cite{Heckman:2012jh,Franco:2012wv}).
Combinatorially, we notice that all of the superpotentials $W$ we have encountered in our examples have the property that each variable (field) appears exactly twice, with opposite sign.
This ensures that each F-term is in the form of ``monomial = monomial'' (so-called toric condition) so that the Jacobian is a so-called binomial ideal, which is one of the equivalent definitions of a toric variety \cite{sturmfels}.

In fact, we can repackage the quiver data, in these examples, into a more compact form: whereas $\cQ$ is a finite graph together with a superpotential $W$ dictating the relations on the path algebra (thus in general these are two independent pieces of information), when $\cM(\cQ)$ is a toric Calabi-Yau variety, however, we can use the following algorithm to recast them into a single bipartite graph on the doubly periodic plane:
\paragraph{Algorithm $\cQ \leftrightarrow \cB$}
\begin{enumerate}
\item Consider a monomial term in $W$, if it comes with a plus (respectively minus) sign, draw a black (respectively white) node \footnote{The choice of colour is, of course, purely by convention.}, write all the variables in the monomial clockwise (respectively counter-clockwise) as edges around the node;
\item Connect all the black/white nodes by these edges.
  Because $W$ has the property that each variable appears exactly twice with opposite sign, we end up with a bipartite graph $\cB$;
\item Because the condition \eqref{T2} is nothing but the Euler relation for a torus, $\cB$ is actually a bipartite graph on $T^2$, or, equivalently, a tiling of the doubly periodic plane;
\item Each node in $\cB$ is a term in $W$ and thus a GIO, and edge is perpendicular to an arrow in $\cQ_1$ obeying orientation and each face in $\cQ$ corresponds to a node in $\cQ_0$; in other words, $\cB$ is a dual graph of $\cQ$;
\item  In particular, being the dual, $\cB$ has $N_2/2$ pairs of black/white nodes, $N_1$ edges and $N_0$ inequivalent (polygonal) faces. 
\end{enumerate}

To illustrate the above procedure, our example from \eqref{c3}, whose VMS is $\IC^3$ yields
\begin{equation}\label{c3dimer}
\begin{array}{ccccc}
\begin{array}{c}
\begin{array}{c}\includegraphics[trim=0mm 0mm 0mm 0mm, clip, width=1in]{c3quiver}\end{array}
\\
W = \tr(XYZ - XZY)
\end{array}
& \longrightarrow
\begin{array}{c}\includegraphics[trim=0mm 0mm 0mm 0mm, clip, width=1in]{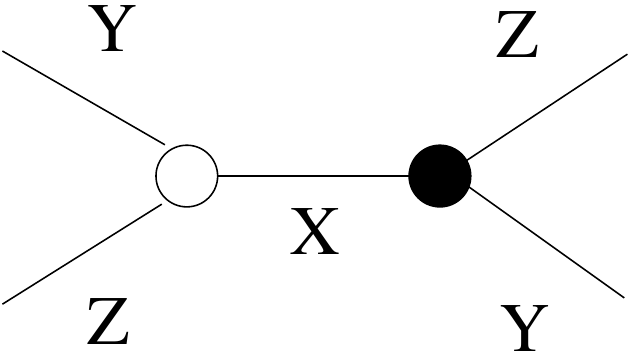}\end{array}
& \longrightarrow
& \begin{array}{c}\includegraphics[trim=0mm 0mm 0mm 0mm, clip, width=1.4in]{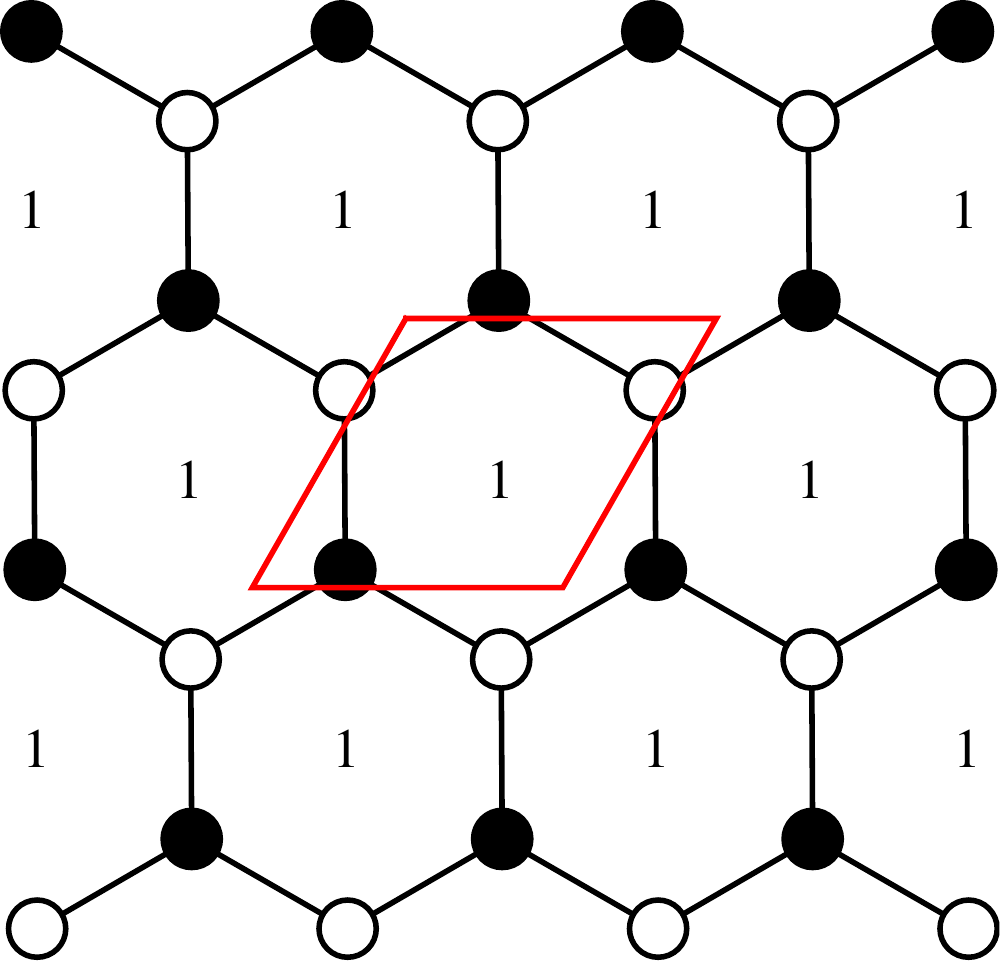}\end{array}
\end{array}
\end{equation}
We have marked the fundamental domain of $T^2$ with the red parallelogram; therein, there is 1 pair of black-white nodes, each of valency 3, corresponding respectively to the $+XYZ$ and $-XZY$ terms in $W$.
The edge together give the honeycomb tiling of the doubly periodic plane, with a single inequivalent face which is a hexagon marked ``1''.

We give another example, the 2-noded quiver with a quartic $W$ is transcribed as
a square tiling of $T^2$, whose fundamental domain consists of 2 pairs of black-white nodes, each of valency 4, corresponding to the 4 terms in $W$ when expanding out the Levi-Civita symbols.
There are 2 inequivalent squares, marked ``1'' and ``2'', corresponding to the 2 nodes in the quiver.
\begin{equation}
\begin{array}{ccc}
  \begin{array}{c}
    \begin{array}{c}\includegraphics[trim=100mm 20mm 50mm 200mm, clip, width=1.5in]{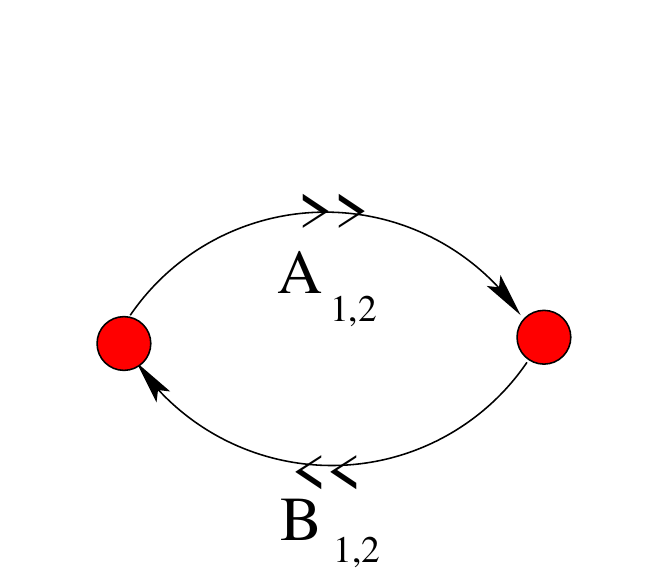}\end{array}
    \\
    W = \tr(\epsilon_{il} \epsilon_{jk} A_i B_j A_l B_k)\\
    \qquad = \tr(A_1B_1A_2B_2 - A_1B_2A_2B_1 + A_1B_2A_2B_1 - A_2B_1A_1B_1) 
\end{array}
& \longrightarrow
& \begin{array}{c}\includegraphics[trim=0mm 20mm 130mm 0mm, clip, width=1.3in]{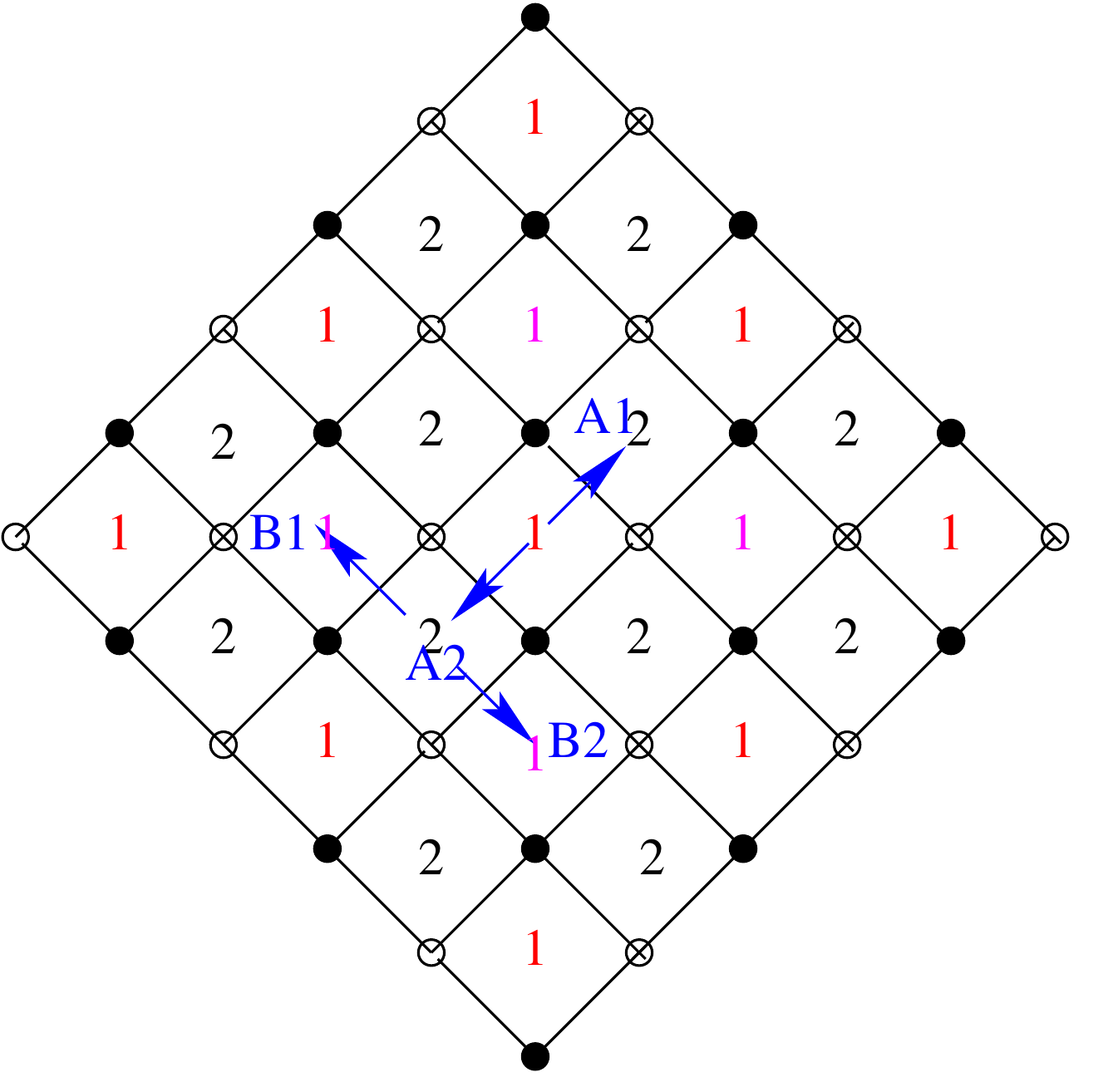}\end{array}
\end{array}
\end{equation}
The VMS $\cM$ here turns out to the be quadric cone $\cC$ to which we alluded earlier and whose toric diagram is given in part (b) of \eqref{toricdiag}.

Over the years, these initial observations by physicists \cite{Feng:2000mi,Feng:2001bn,Feng:2005gw,Franco:2005rj,Franco:2005sm,Hanany:2005ve} have been formalized into the following result by a host of mathematicians \cite{AU,beil,BKM,BCQ,rafCY,Broom,davison,DWZ,Ginz,ishii2010note,LM,MR,balazs}
\begin{theorem}
  Let $\cQ$ be a quiver graph-dual to a bipartite graph $\cB$ drawn on $T^2$ as outlined in the algorithm above, then the (coherent component of) $\cM(\cQ)$ is an affine toric Calabi-Yau variety of dimension three.
\end{theorem}
In general, we remark that for a bipartite graph drawn on a (compact) genus $g$ Riemann surface without boundaries (so physically this means we are not considering flavour nodes corresponding to various matter in the  fundamental representation), the VMS is a toric Calabi-Yau variety of dimension $2g+1$.

Having established this reassuring correspondence between bipartite graphs and Calabi-Yau varieties, one can then use the rich structures, combinatorial from the former and geometrical, the latter, to arrive a plenitude of consequences.
We only mention a few major corollaries here.

\paragraph{Kasteleyn Matrix: }
The most important object in the study of bipartite graphs, either in combinatorics or in the physics of condensed matter theory (wherein our graphs are referred to as {\it ``dimers''}), is the Kasteleyn matrix \cite{kas}, which computes the number of ways of completely pairing black/white nodes, known as {\it perfect matchings} \footnote{
  In fact, with the definition of a height function, the determinant of the Kasteleyn matrix actually encodes the number of perfect matchings in the dimer for given height change; this in turn corresponds exactly to the multiplicities of GLSM fields in the gauge theory \cite{Feng:2000mi,Franco:2005rj,Franco:2005sm}, though this beautiful side of the story we will not presently need.
  }.
Given a bipartite graph $\cB$ on $T^2$, the Kasteleyn matrix is constructed as a weighted adjacency matrix in the following manner:
\begin{enumerate}
\item There is a natural orientation induced from the bipartite colouring: choose, say, ``black'' to ``white'';
\item Weight every edge by $e_i = \pm 1$ so that for {\it each} face, bounded by say $E$ edges $e_{i=1,\ldots, E}$, we have Sign$(\prod e_i) = +1$ if $E = 2 (\bmod 4)$ and $-1$ if $E = 0 (\bmod 4)$; this is always doable;
\item Draw 2 paths in the fundamental domain of $T^2$ as $\gamma_w$ and $\gamma_z$ corresponding to the 2 independent homology cycles of the torus;
\item Set up an $(N_2/2) \times (N_2/2)$ matrix $K_{ij}$ whose rows, say, index the white nodes and whose columns, black nodes.
  If white node $i$ is not connected to black node $j$ by an edge, then $K_{ij} = 0$. If there is an edge and it crosses $\gamma_z$, then we introduce a factor of $z$ or $1/z$ according to our chosen orientation (likewise, $w$ or $1/w$ for $\gamma_w$). Finally, multiply this entry in $K_{ij}$ by our chosen sign of the edge itself in step 2;
\item If there are multiple edges between the nodes, we add the contributions.
\end{enumerate}

As an illustration, we use the following example, taken from Fig.~4 of \cite{Franco:2005rj}:
\[
\begin{array}{c}
    \includegraphics[trim= 0mm 0mm 0mm 0mm,clip,width=6in]
                    {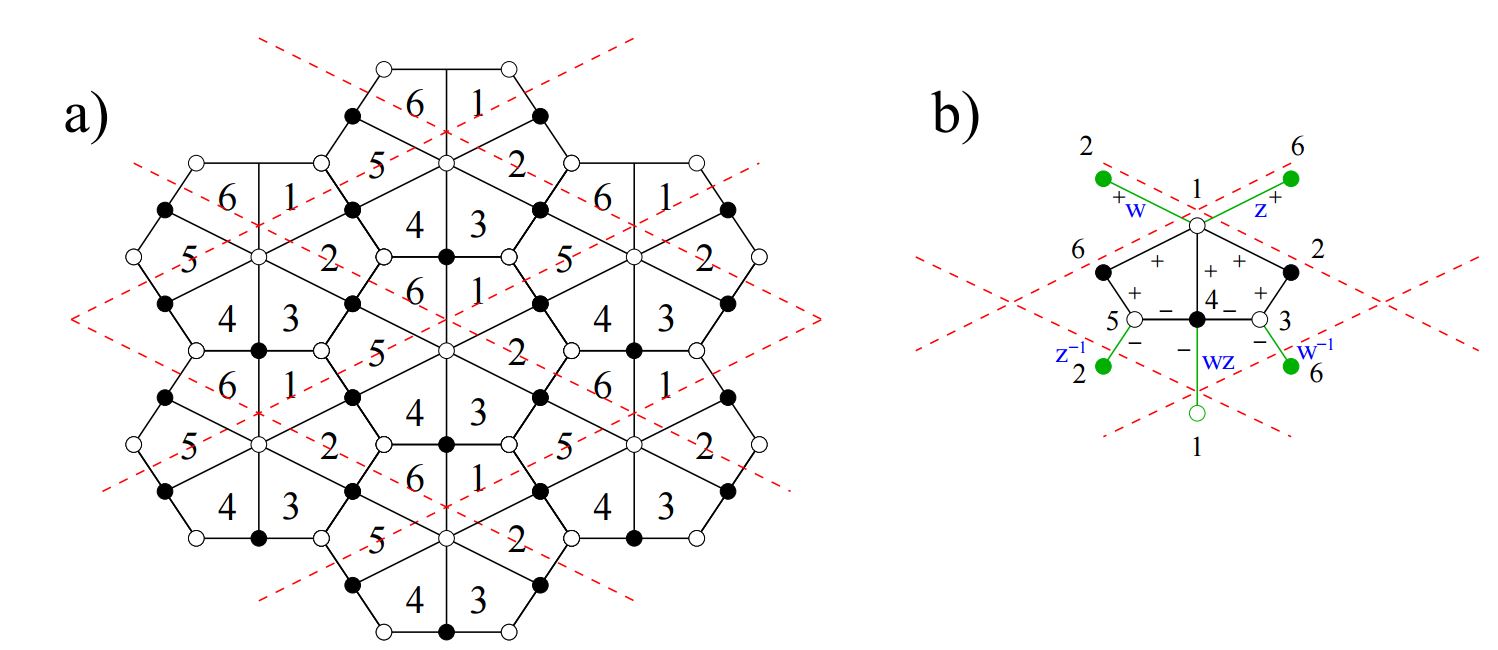}
\end{array}
\]
Suppose we are given the bipartite graph in part (a), which has $N_2/2 = 3$ pairs of white-black nodes (marked as 1,3,5 and 2,4,6 respectively), $N_1 = 12$ edges and $N_0 = 6$ faces, then, zooming into the fundamental domain in part (b), with $\gamma_z$ and $\gamma_w$ indicated and the factors of $z^{\pm1}$, $w^{\pm1}$, as well as the sign assignment of the edges explicitly marked, one can check that we arrive at the $K(z,w)$ below, which is $3 \times 3$ and has 12 total monomial terms:
\begin{equation}
  K(z,w) = \left(
  \begin{array}{c|ccc}
      & 2 & 4 & 6 \\ \hline
    1 & 1+w & 1-zw & 1+z \\
    3 & 1 & -1 & -w^{-1} \\
    5 & -z^{-1} & -1 & 1
  \end{array}
  \right) \
  \det(K(z,w)) = w^{-1}z^{-1} - z^{-1} - w^{-1} - 6 - w - z + wz \ .
\end{equation}

In the above, we have presented the determinant of $K(z,w)$ for reference but the vigilant reader would perhaps reflect upon Definition \ref{newton}.
The determinant here is the Newton polynomial for the toric diagram for $dP_3$ (cf.~\S\ref{s:CY3}) and indeed the quiver in part (a) is that associated with this geometry.
In general, we have \cite{Franco:2005rj,Feng:2005gw}
\begin{theorem}
  The determinant of the Kasteleyn matrix of $\cB$ is the Newton polynomial of the toric diagram for the affine Calabi-Yau threefold $\cM(\cQ)$,
  \[
  \det(K_{\cB}(z,w)) = P_{\cD(\cM(\cQ))}(z,w) \ .
  \]
\end{theorem}
Again, $\cB$ and $\cQ$ are dual graphs outlined in our algorithm above.
The coefficients in $\det(K(z,w))$ are, by construction, integers and these count the perfect matchings/GLSM field multiplicities; they are indeed special points in the moduli space of the defining equation for the mirror of $\cM$.

It is amusing to reflect upon the history.
Kasteleyn obtained his matrix almost a decade before the discovery of toric varieties and yet his determinant is none other than the Newton polynomial.
One might wonder whether he ever thought of $\det(K_{\cB}(z,w))$ as encoding a lattice polygon and if so how it would have befuddled him.

\paragraph{Quivering Amoebae: } 
Given an affine variety, especially a Riemann surface, there are two well-known {\it real} projections as algebraic sets.
These are complementary, and essentially project to the modulus and phase parts of the complex variables \cite{kos,uy,Feng:2005gw}:
\begin{definition}
  For a Riemann surface $\Sigma: P(z,w) = 0$, there are two projections
  \begin{enumerate}
  \item Amoeba $Am(\Sigma) := \big(\log(|z|), \log(|w|)\big) \subset \IR^2$;
  \item Coamoeba $CoAm(\Sigma) :=  \big(\arg(z), \arg(w)\big) \subset T^2$ \ .
  \end{enumerate}
\end{definition}
In the above, the coamoeba lives on a torus because the argument is defined on $[0,2\pi)$, moreover, in \cite{Feng:2005gw}, it was called the ``alga''.
  One can also see readily that the above definition generalizes to any complex variety, but we will only need it for the bivariate case.

  With these definitions, we can complete the cycle of relations between the quiver, the dimer, and the toric variety:
  \begin{theorem}
    Given $\cQ$ and its dual graph $\cB \subset T^2$, whose CY3 VMS $\cM(\cQ)$ has toric diagram $\cD$, we have that
    \begin{enumerate}
    \item \cite{Feng:2005gw} The deformation retract (spine) of $Am(P_{\cD}(z,w))$ is the graph dual to $\cD$;
    \item \cite{uy} The deformation retract of $CoAm(P_{\cD}(z,w))$ is $\cB$.
    \end{enumerate}
  \end{theorem}
  Here, the Newton polynomial $P_{\cD}(z,w)$, by foregoing discussions, is obtained from the Kasteleyn matrix, and also determines the complex equations of the local mirror to $\cM$.
  Furthermore, the dual graph to the planar toric diagram is called (p,q)-web in physics.

\paragraph{The Plethystic Programme: }
In the course of investigations of the quivers as gauge theories, and motivated by enumeration of certain GIOs, wherein a programme of utilizing the Hilbert series of $\cM$ was established, some intriguing properties of the Hilbert series of algebraic varieties in general were observed \cite{Benvenuti:2006qr}.
We will present some key points as observations since a rigourous treatment is yet to be administered.

We begin with the definition of plethystics
\begin{definition}
  Given a smooth function $f(t)$, its plethystic exponential is the formal series
  \[
  PE[f(t)] = \exp\left( \sum_{n=1}^\infty \frac{f(t^n) -
  f(0)}{n} \right) \ .
  \]
\end{definition}
Assuming reasonable regions of convergence, one can readily check from this definition by direct manipulation of the series expansions, that
\begin{lemma}
  If $f(t)$ has Taylor series $f(t) = \sum\limits_{n=0}^\infty a_n t^n$, then
  \begin{itemize}
  \item There is an Euler-type product formula
    $g(t) = PE[f(t)] = \prod\limits_{n=1}^\infty (1-t^n)^{-a_n}$;
  \item There is explicit inverse (called plethystic logarithm) such that
    $f(t) = PE^{-1}(g(t)) = \sum\limits_{k=1}^\infty \frac{\mu(k)}{k} \log (g(t^k))$
  \end{itemize}
\end{lemma}
In the above, $\mu(k)$ is the standard M\"obius function for $k \in \IZ_+$ which is 0 if $k$ has repeated prime factors and is $(-1)^n$ if $k$ factorizes into $n$ distinct primes (also with the convention that $\mu(1)=1$.

Thus prepared, we can specialize to $f(t)$ being a Hilbert series, whose Taylor coefficients are by definition non-negative integers.
\begin{observation}\label{plog}
  Given Hilbert series $H(t; \cM)$ of an algebraic variety $\cM$, the plethystic logarithm is of the form
  \[
  PE^{-1}[H(t; \cM)] = b_1 t + b_2 t^2 + b_3 t^3 + \ldots
  \]
  where all $b_n \in \IZ$ and a positive $b_n$ corresponds to a generator in coordinate ring of $\cM$ and a negative $b_n$, a relation.
  In particular, if $\cM$ is complete intersection, then $PE^{-1}[H(t; \cM)]$ is a finite polynomial.
\end{observation}
For example, our quadric hypersurface $\cC$ in $\IC^4$ has Hilbert series
$H(t; \cC) = (1 - t^2) / (1 - t)^4$ and one checks that $PE^{-1}[H(t; \cC)] = 4t - t^2$, meaning that there are 4 generators in degree 1 (corresponding to the 4 complex coordinately of $\IC^4$) with 1 relation at degree 2 (meaning that $\cC$ is a single quadratic hypersurface, and hence complete intersection, in $\IC^4$).
While for this example, the result trivially follows from the very construction of the Hilbert series, whether $PE^{-1}$ only terminates for complete intersections and how one might disentangle positive and negative contributions to each $b_n$ is, as far as we are aware, not in general known.

\newpage

\subsection{Cluster Mutation and Seiberg Duality}
Perhaps the most extraordinary result, for which we reserve this separate subsection, is the preservation of $\cM(\cQ)$ under certain transformations of $\cQ$, motivated independently from the physics and from the mathematics.
In a supersymmetric gauge theory, it is a classic result \cite{Seiberg:1994pq} that an $SU(N_c)$ theory with $N_f$ flavours $(Q,Q')$ is quantum-field-theoretically dual to an $SU(N_f - N_c)$ theory with $N_f$ flavours $(q,q')$ and superpotential $W = Mqq'$, for $\frac32 N_c \leq N_f \leq 3 N_c$; this is called {\it Seiberg Duality}.
Now, though this sentence would mean very little to the mathematician, it was realized in \cite{Feng:2000mi,Feng:2001bn,Beasley:2001zp,Cachazo:2001sg}, that given our dictionary at the very outset of this review, this duality is actually a very particular transformation on the quiver.
In particular, Seiberg duality corresponds to
\begin{proposition}\label{def:seiberg}
  Seiberg duality, on a given node with $N_c$ as its label (component of dimension vector) with $N_f$ incoming and $N_f$ outgoing arrows, is the following transformation  on the quiver, where the label is changed, arrows are reversed, and a superpotential is generated:
  \[
  \begin{array}{c}
    \includegraphics[trim=0mm 0mm 0mm 0mm, clip, width=5in]{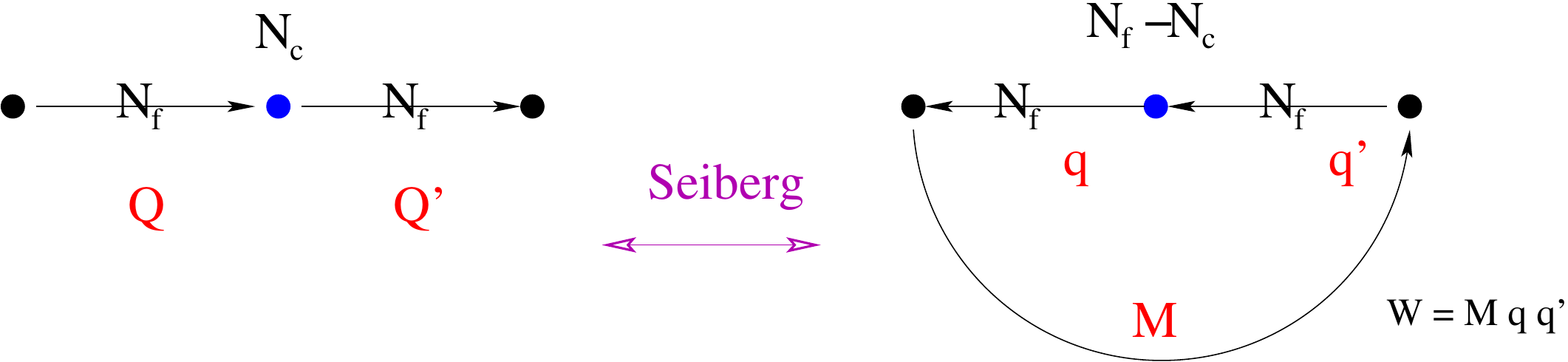}
  \end{array}
  \]
\end{proposition}
We can take this to be our working definition of Seiberg duality and we can easily check that it is an involution on the space of quivers.
Now, recalling that the dimension vector for our quivers are all $(1,1, \ldots, 1)$ and we wish to stay within this realm for convenience of discussion, we will thus focus on Seiberg duality on nodes which have 2 incoming and 2 outgoing arrows (which has been variously dubbed toric duality \cite{Feng:2000mi} or toric Seiberg duality \cite{Hanany:2011ra}).

\paragraph{Remark: }
The informed reader would recognize the above transformation rule as nothing but the {\bf cluster mutation} rules of \cite{FZ} (cf.~\cite{DWZ,Mar}), coming from an entirely different motivation.
This is quite remarkable. In 2000, the physicists and the mathematician were completely unaware of each others' work and yet almost contemporaneously, they hit upon these rather peculiar transformation rules on a quiver $\cQ$.
Only years later did the dialogue between the two disparate communities working on the same subject begin over a series of workshops.

From the dual point of view of the bipartite graph $\cB$, one can readily translate this cluster mutation rule \cite{Franco:2005rj,Franco:2005sm}:
\begin{theorem}\label{seiberg}
  Seiberg duality/quiver mutation on $\cQ$ corresponds to the following transformation on the dimer model/brane tiling $\cB$
  \[
  \begin{array}{c}
    \includegraphics[trim=0mm 0mm 0mm 0mm, clip, width=4.5in]{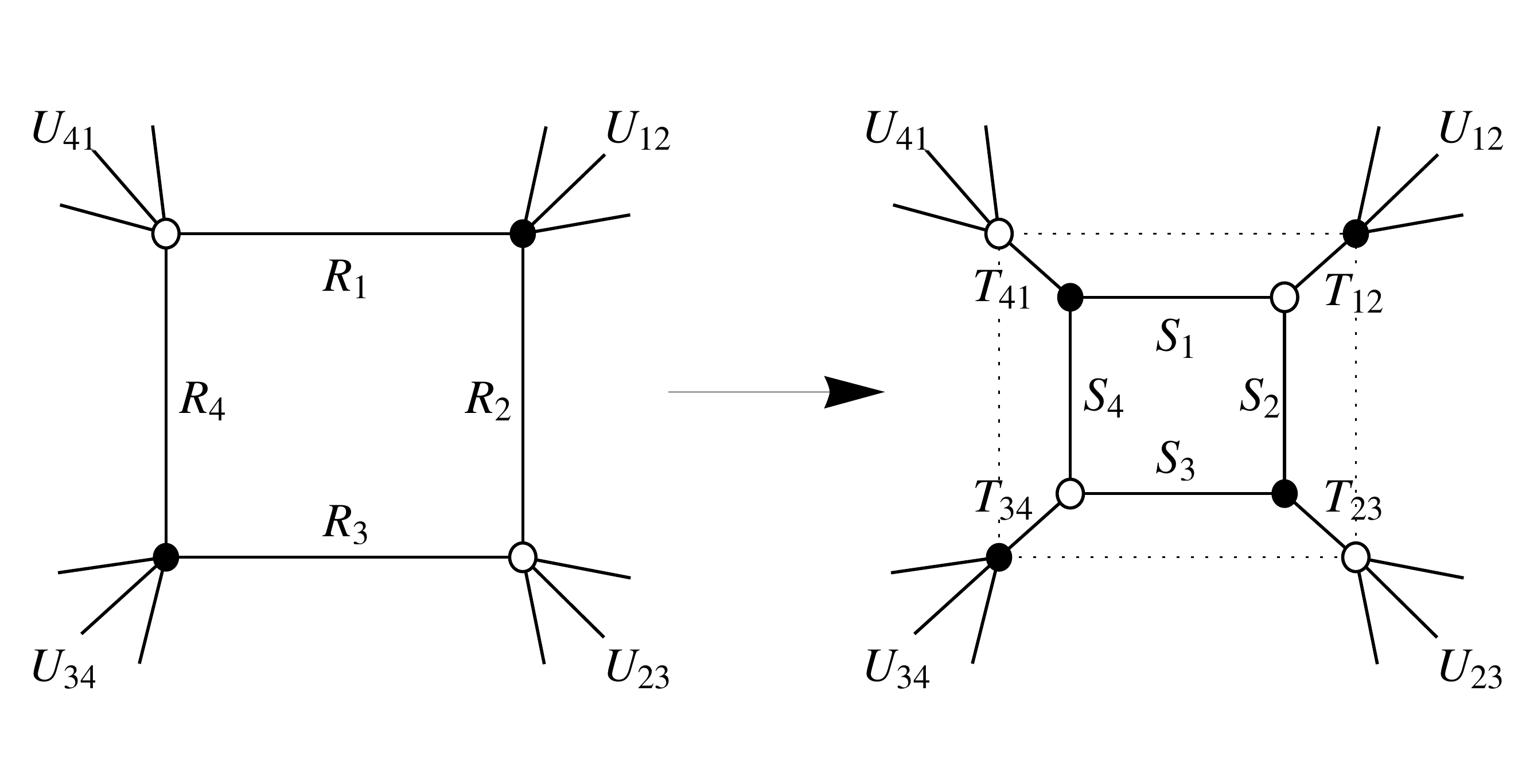}
  \end{array}
  \]
  which, jocundly, has been called ``urban renewal''.
\end{theorem}
We see that under this renewal process, $\cB$ remains bipartite, new edges are created (corresponding to the new arrows generated in the quiver under mutation) and new vertices (corresponding to additional superpotential terms) emerge.

Both Seiberg duality, being a quantum field theory duality, and cluster mutation, being generator of a cluster algebra, have a plethora of important consequences.
For our present purposes on studying the VMS, a pertinent result is
\cite{Feng:2000mi,Franco:2005rj,Franco:2005sm}
\begin{corollary}\label{equalM}
  For quivers $\cQ$ and $\cQ'$ related by Seiberg duality/cluster mutation (or, equivalently, dimers $\cB$ and $\cB'$ related by urban renewal), we have that
  \[
  \cM(\cQ) \simeq \cM(\cQ')
  \]
  as affine varieties.
\end{corollary}
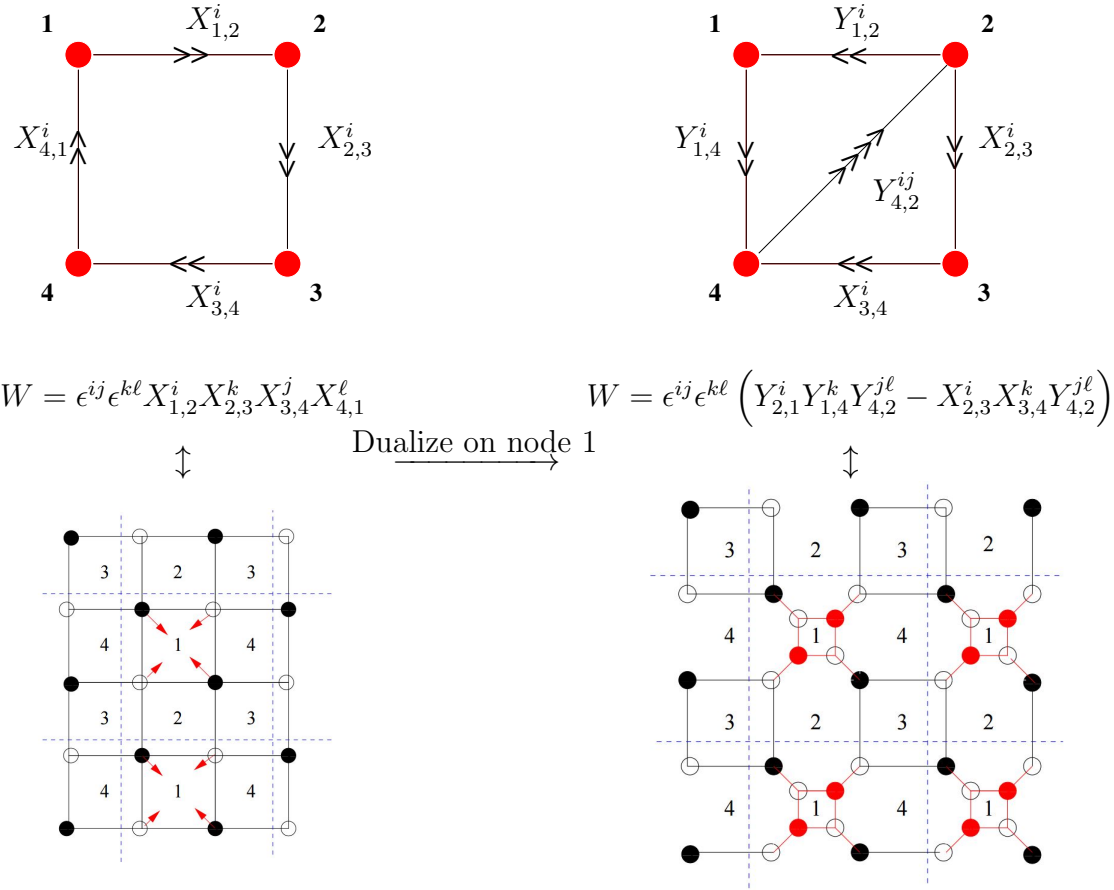
\begin{figure}[t!]
\begin{picture}(500,350)(0,-150)
  \put(75,165){$X_{1,2}^i$}
  \put(75,60){$X_{3,4}^i$}
  \put(10,120){$X_{4,1}^i$}
  \put(125,120){$X_{2,3}^i$}
  \put(320,165){$Y_{1,2}^i$}
  \put(320,60){$X_{3,4}^i$}
  \put(260,120){$Y_{1,4}^i$}
  \put(375,120){$X_{2,3}^i$}
  \put(335,100){$Y_{4,2}^{ij}$}
  $
  \begin{array}{cccc}
  \begin{array}{c} \includegraphics[trim= 0mm 0mm 0mm 0mm,clip,width=1.5in]{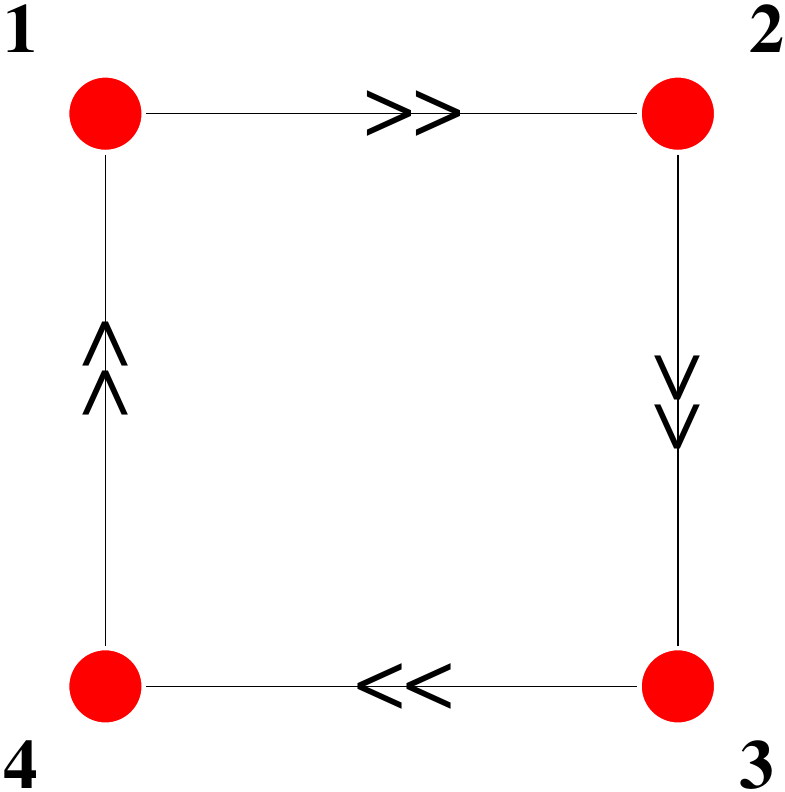}\end{array}
  &
  &
  \begin{array}{c} \includegraphics[trim= 0mm 0mm 0mm 0mm,clip,width=1.5in]{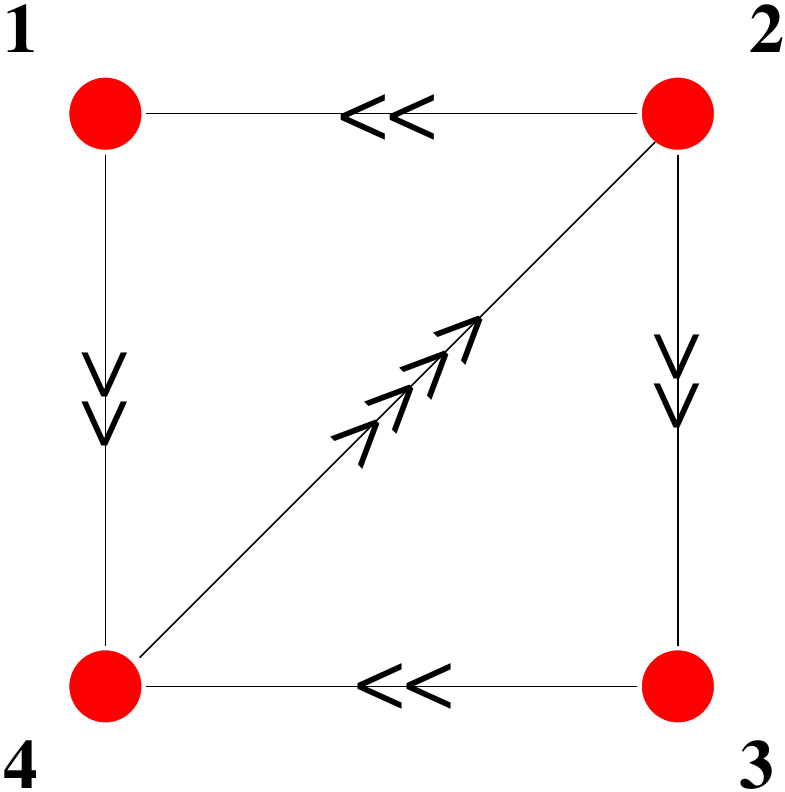}\end{array} \\
  &&&\\
  W = \epsilon^{ij}\epsilon^{k\ell} X_{1,2}^i X_{2,3}^k X_{3,4}^j X_{4,1}^{\ell}
  &
  &
  W = \epsilon^{ij}\epsilon^{k\ell} \left(
  Y_{2,1}^i Y_{1,4}^k Y_{4,2}^{j\ell} - X_{2,3}^i X_{3,4}^k Y_{4,2}^{j\ell}
  \right)
  \\
  \updownarrow &  
  \xrightarrow{\makebox[2cm]{\mbox{Dualize on node 1}}}  &
  \updownarrow
  \\
  \begin{array}{c} \includegraphics[trim= 0mm 0mm 0mm 0mm,clip,width=1.6in]{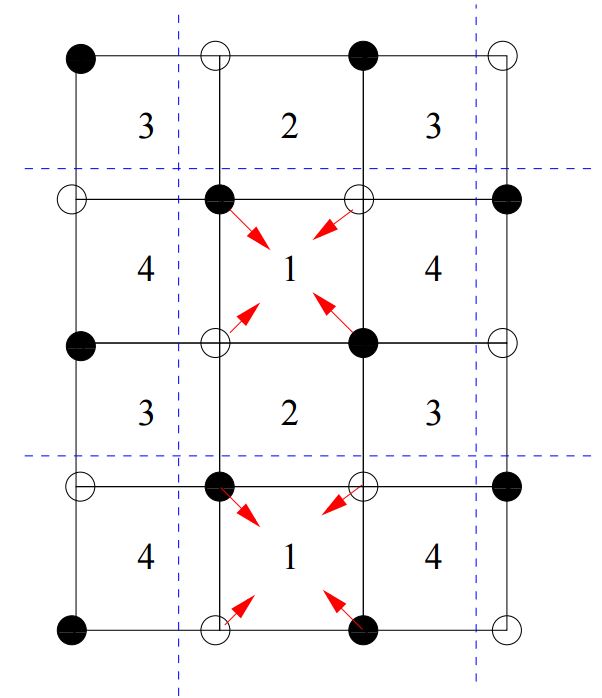}\end{array}
  &
  &
  \begin{array}{c} \includegraphics[trim= 0mm 0mm 0mm 0mm,clip,width=2.3in]{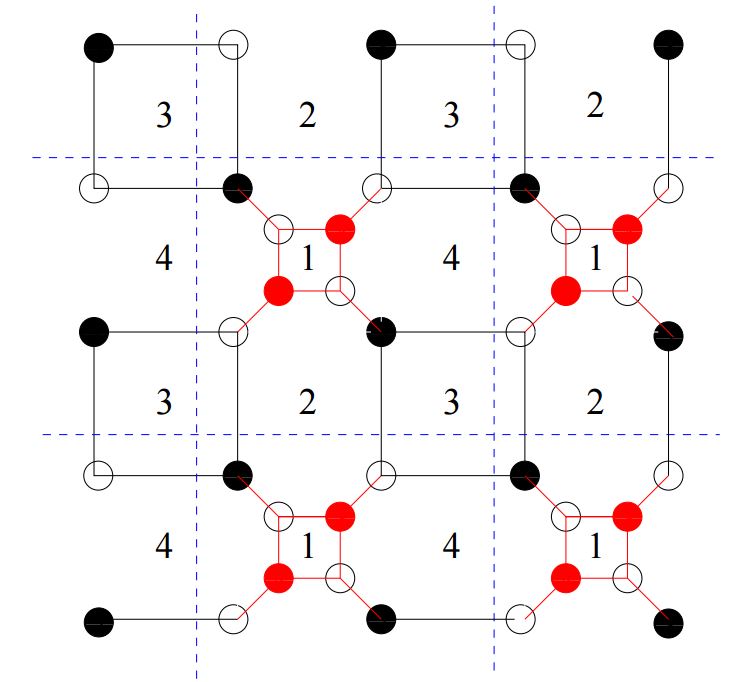}\end{array} \\
  \end{array}
$
\end{picture}
\caption{{\sf
  The Seiberg dual/cluster transformation pair of quivers (equivalently, urban renewal on the pair of dimers on $T^2$) whose VMS are both $\IF_0$, the cone over $\IP^1 \times \IP^1$.}
  \label{f:F0}}
\end{figure}

There is by now a host of examples of (toric) Seiberg duals, the most archetypical of which is the pair for $\IF_0$ (which we recall from \S\ref{s:CY3} is the CY3 cone over $\IP^1 \times \IP^1$), which is presented in Figure \ref{f:F0}.
Suppose we have the quiver on the left with 8 arrows which we can label as $X_{a,a+1}^{i=1,2}$ with $a = 1,2,3,4$ and defined modulo 4.
The quiver, superpotential and the dual dimer are given to the left.
Now, let us dualize on node 1 according to the rule in Proposition \ref{def:seiberg}, whereupon we arrive at the data to the right, with arrows $Y_{1,2}^i$ and $Y_{1,4}^i$ the reverse of $X_{1,2}^i$ and $X_{4,1}^i$ respectively and a new set of 4 arrows $Y_{4,2}^{ij} \sim X_{4,1}^iX_{1,2}^j$ with $i,j=1,2$ being generated.
Additional superpotential terms in the form of $Y_{2,1}^i Y_{1,4}^k Y_{4,2}^{j\ell}$ is also generated.
From the perspective of the dual dimer, the square tiling of $T^2$ is transformed into the square-octagon tiling.
Crucially, one can check, using the method \footnote{
  In fact, as one can imagine, for toric VMS, there is a much more efficient and purely combinatorial method of obtaining $\cM$. 
} of Proposition \ref{vms}, that $\cM$ for both quivers is the affine Calabi-Yau cone over $\IP^1 \times \IP^1$.

\section{Dessins d'Enfants: the Shape of $T^2$}\label{s:dessin}
Whilst in the foregoing we have discussed much about the combinatorics of the bipartite graph, the subject of {\it shape} and {\it size} of the $T^2$ which our objects tile have not arisen.
Topologically, that we are dealing with $U(1)^n$ quivers, combined with the ``toric condition'' that each arrow appears in $W$ exactly twice with opposite sign, together with the ``conformality'' condition \eqref{T2}, compel us to be drawing bipartite graphs on $T^2$.
Though it was further realized \cite{Feng:2005gw} that the $T^2$ is part of the $T^3$ fibration in local mirror symmetry \cite{Strominger:1996it}, and that the size of the $T^2$ is, as required in string compactification, on the microscopic scale of Planck length, its shape requires further input from physics.

\subsection{Isoradial Embeddings}
In supersymmetric quantum field theories, the space-time fields carry a charge from the supersymmetry representations, called {\em R-Charge}; in the language of quivers/tilings, we can take the following as a working definition \footnote{
  The first condition comes from the fact that all superpotential terms must carry R-charge 2 in order to be well-defined in superspace and the second condition is a consequence of conformality so that the NSVZ beta-function vanishes for each gauge group. 
  }
\begin{definition}\label{Rcharge}
  To each edge in the bipartite graph $\cB$, one can assign a number $R_i \in \IR_+$, called R-charge, such that
  \begin{enumerate}
  \item For each node, $\sum\limits_{\substack{i \in \text{ edges} \\ \text{around node}}} R_i = 2$;
  \item For each face, $\sum\limits_{\substack{i \in \text{ edges}\\  \text{bounding face}}} (1-R_i) = 2$.
  \end{enumerate}
\end{definition}
In fact, we can readily verify (by summing the first equation over the nodes, and the second, over the faces) that this pair of conditions suffice to imply the topological condition \eqref{T2}.

\begin{figure}[h!]
  (a)
  \begin{tabular}{c}\includegraphics[trim= 0mm 0mm 0mm 0mm,clip,width=2.5in]{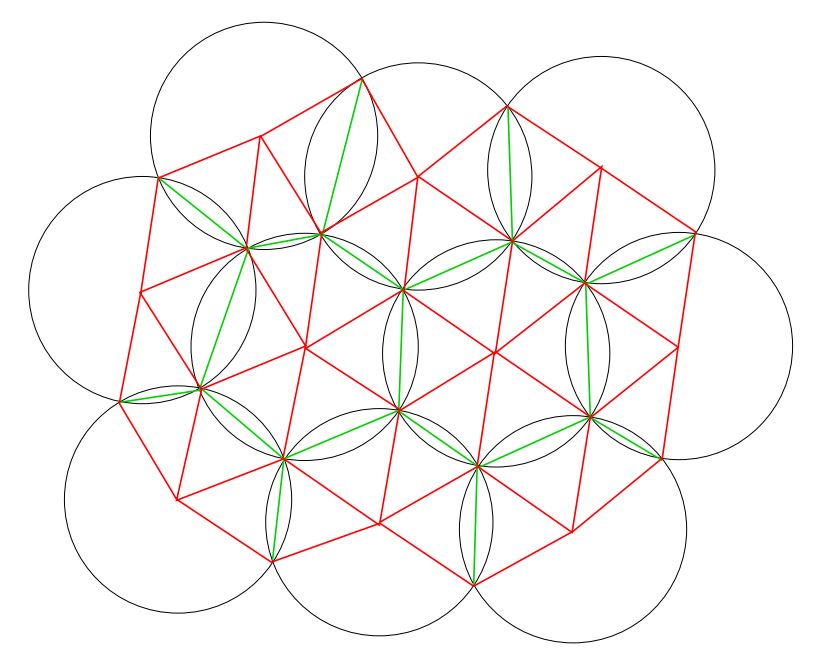}\end{tabular}
  (b)
  \begin{tabular}{c}\includegraphics[trim= 0mm 0mm 0mm 0mm,clip,width=2.5in]{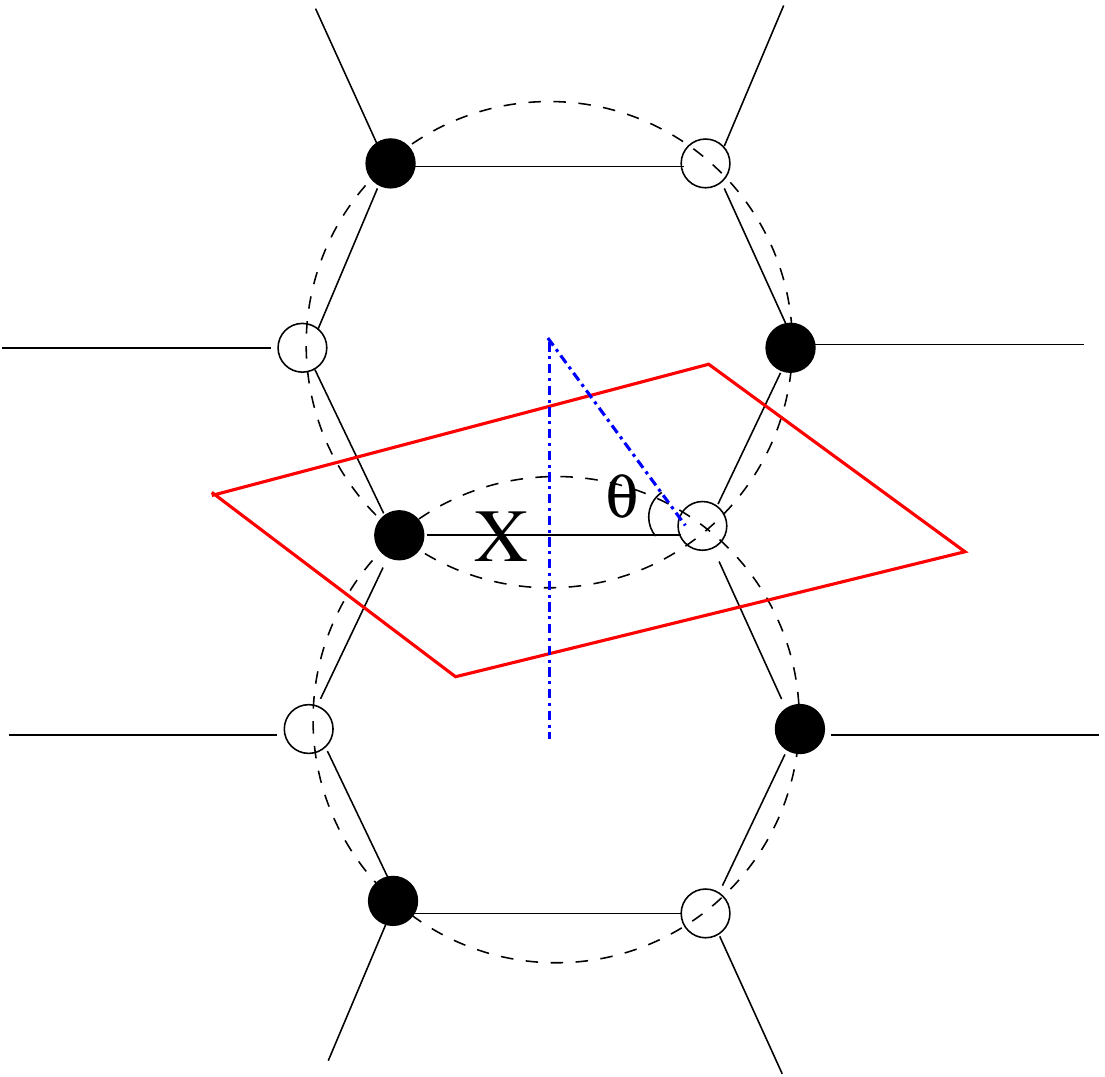}\end{tabular}
\caption{{\sf
    (a) In an isoradial embedding, the edges of the dimer, together with the radius of the circles, form various rhombi.
    (b) Setting $R_i = \frac{2 \theta_i}{\pi}$ satisfies the constraints in Definition \ref{Rcharge}.
    }
  \label{f:isoradial}}
\end{figure}

This assignment thus prescribes a natural (from the point of view of physics) length to the edges in $\cB$ and there is a simple way to satisfy both constraints \cite{Hanany:2005ss} (which can be seen by multiplying the 2 condition by $\pi$ and elementary geometry):
\begin{lemma}
  Let all black/white nodes be drawn on circles of equal radius (cf.~part (a) of Figure \ref{f:isoradial}) in a so-called isoradial embedding. In this way, every quadrilateral formed by a pair of black/white nodes and the centres of adjacent circles is a rhombus.
  Let $2\theta_i$ be the angle subtended by the isosceles formed by the adjacent centres and a node, hence the line through the two centres is the perpendicular bisector of an edge $X_i$ (cf.~part (b) of said Figure).
  Then
  \begin{equation}\label{Rtheta}
    R_i = \frac{2 \theta_i}{\pi}
  \end{equation}
  satisfies the conditions in Definition \ref{Rcharge}.
\end{lemma}

We still have one last degree of freedom for the nodes in $\cB$: they can still move along the circumference of the isoradial circles.
Again, there is a final physics constraint \footnote{
  The origin of this function is the trial a-function $a := \frac{3}{32}( 3 \tr R^3 - \tr R)$ but for us $\tr R = 0$.}
called a-maximization \cite{Intriligator:2003jj} which dictates that when computing the VMS (i.e., in the infra-red), one has \cite{Martelli:2006yb,Butti:2005vn,Kato:2006vx,Hanany:2011bs} 
\begin{theorem}\label{amax}
  In the VMS $\cM(\cQ)$, the R-charge assignment above needs to maximize the real function
  \begin{equation}\label{aR}
  a(R_i) = \sum\limits_{\text{all edges}} (R_i - 1)^3 \ ;
  \end{equation}
  and the maximum is unique.
\end{theorem}
In fact, this function corresponds to the normalized volume of the base Sasaki-Einstein manifold to $\cM(\cQ)$ as shown in \cite{Martelli:2006yb} and is related to the Futaki-Donaldson invariant \cite{Collins:2016icw}.

In summary, the isoradial embedding of our bipartite graph $\cB$ on $T^2$, such that the R-charge assignments \eqref{Rtheta} maximize \eqref{aR}, fixes the {\em shape} of $T^2$ completely.
The torus is thus endowed with complex structure $\tau$ with Klein invariant $j(\tau)$, and on which $\cB$ is a rigidly embedded graph. 
Its {\em size}, on the other hand, is not specified -- we could have chosen the radii of circles in Figure \ref{f:isoradial} arbitrarily; this corresponds to K\"ahler structure which ultimately determines the complex structure of the mirror to $\cM$.
As mentioned earlier, this size is taken to be microscopically small.

Let us illustrate the above with our running examples.
For $\IC^3$, there are three inequivalent edges, thus we need to solve for $R_{i,1,2,3}$ in order to find
\begin{eqnarray}
  \nn
  &\text{Max}\left[ a(R_i) = (R_1 - 1)^3 + (R_2 - 1)^3 + (R_3 - 1)^3 \right] \ ,
  \\
  &
  \mbox{ such that }
  R_1 + R_2 + R_3 = 2 \ , \quad
  2\big( (1 - R_1) + (1 - R_2) + (1 - R_3) \big) = 2 \ .
\end{eqnarray}
Note that the second constraints has a factor of 2 because the 3 edges form a hexagon and thus each appears twice in the face.
We readily solve, using Lagrange multipliers, to arrive at
$R_1=R_2=R_3 = \frac23$.
Drawing these as angles in $\cB$ renders
\begin{equation}\label{aMaxC3}
  \theta_1 = \theta_2 = \theta_3 = \frac{\pi}{3}
  \Rightarrow
  \tau = -\frac12 + \frac{\sqrt{3}}{2} i
  \Rightarrow
  j(\tau) = 0 \ .
\end{equation}

Next, let us try the two Seiberg/cluster dual pairs of $\IF_0$ in Figure \ref{f:F0}, which we call $\IF_0(I)$ and $\IF_0(II)$.
In $\IF_0(I)$, we have 8 edges and 4 faces, thus we need to find $R_{i=1,\ldots,8}$ satisfying 4 each of the two constraints in Definition \ref{Rcharge} - the solution space of which is actually 4 dimensional, say parametrized by $R_{1,2,3,4}$ - and maximizing the cubic $2 \sum\limits_{i=1}^4 (R_i - 1)^3$.
We will not present the tedium of the details, but only the result:
\begin{equation}\label{F0IR}
  R_i = \frac12
  \Rightarrow
  \theta_i = \frac{\pi}{4}
  \Rightarrow
  \tau = i
  \Rightarrow
  j(\tau) = 1 \ .
\end{equation}
Similarly, in $\IF_0(II)$, we have 12 edges and 4 faces, thus we need to solve for $R_{i=1,\ldots, 12}$ which must first obey the 8 + 4 conditions of Definition \ref{Rcharge} - which again are not all independent and affords a 9 dimensional space of solutions - and maximizing the cubic in \eqref{aR}, yielding
\begin{equation}\label{F0IIR}
  R_{i=1,\ldots,8} = \frac12, \ R_{i=9,\ldots,12} = 1
  \Rightarrow
  \tau = i
  \Rightarrow
  j(\tau) = 1 \ .
\end{equation}

One immediately observes that $j(\tau)$ is the same for both $F_0(I)$ and $F_0(II)$, and might suspect that the {\em shape} of the $T^2$ to be invariant under Seiberg duality; and this is indeed correct.
Via the rigidification of the $T^2$ as determined by Definition \ref{Rcharge} and maximization of \eqref{aR}, the $T^2$ becomes an elliptic curve with a particular $j$-invariant and the R-charges are all algebraic numbers since we are solving algebraic equations throughout.
The general statement is that (q.v.~\S3 of \cite{Hanany:2011bs} for the demonstration)
\begin{theorem}\label{clusterj}
  Under Seiberg duality/cluster mutation/urban renewal, the following are invariants
  \begin{enumerate}
  \item
    $j(\tau)$ of the $T^2$ which $\cB$ tiles;
  \item
    the transcendence degree of the R-charges over $\IQ$.
  \end{enumerate}
\end{theorem}
In the second part, by transcendence degree we mean the degree of the minimal field extension over $\IQ$ which contains all the R-charges.
We will appreciate this purely number theoretic quantity a little more in the next subsection.

\subsection{Grothendieck's Dessin d'Enfants}
When one infringes upon the subject of bipartite graphs on $T^2$, or embedded graphs on Riemann surfaces in general, one inevitably is led to the subject of Grothendieck's {\em dessins d'enfants}, or children's drawings, one of the most profound pieces of modern mathematics \cite{leila}, residing in the intersections between number theory, geometry and combinatorics.
Indeed, in his own immortal words,
\begin{quote}
This discovery, which is technically so simple, made a very strong impression on me, and it represents a decisive turning point in the course of my reflections, a shift in particular of my centre of interest in mathematics, which suddenly found itself strongly focused. {\em I do not believe that a mathematical fact has ever struck me quite so strongly as this one, nor had a comparable psychological impact.}

This is surely because of the very familiar, non-technical nature of the objects considered, of which any child's drawing scrawled on a bit of paper (at least if the drawing is made without lifting the pencil) gives a perfectly explicit example. To such a dessin we find associated subtle arithmetic invariants, which are completely turned topsy-turvy as soon as we add one more stroke.
\end{quote}

Let us remind ourselves of some key results of dessins, starting with the extraordinary theorem of Belyi \cite{belyi}:
\begin{theorem}
Let $\Sigma$ be a smooth compact Riemann surface, then $\Sigma$ has an algebraic model over $\overline{\IQ}$ IFF there is surjective map $\beta: \Sigma \to \IP^1$ ramified at exactly 3 points.
\end{theorem}
There are several points of note about this theorem.
First, the IFF translates the analytic property of ramification (i.e., in local coordinates the Taylor series of $\beta$ starts from order 2) to the number-theoretical property of the algebraic closure of $\IQ$.
Second, the proof, which uses classical methods essentially known to Riemann, arose as late as 1980.

Using $SL(2;\IC)$ on $\IP^1$, we can take the three branch points to be $(0,1,\infty)$.
Then the theorem states that $\Sigma$, defined as a polynomial in $\IC[x.y]$ (say, as a hyper-elliptic curve $y^2 = p(x)$), has coefficients as algebraic numbers IFF one can find a rational function $f(x,y) : \Sigma \to \IP^1$ whose Taylor series has no linear term at $0,1,\infty$.

Grothendieck's insight was to realize that the Belyi map gives an embedded graph on $\Sigma$ as follows:
\begin{definition}
  Consider $\beta^{-1}(0)$, which is a set of points on $\Sigma$ that can be coloured as black, and likewise $\beta^{-1}(1)$, white.
  The pre-image of any simple curve with endpoints 0 and 1 on $\IP^1$ is a bipartite graph embedded in $\Sigma$ whose valency at a point is given by the ramification index (i.e., order of vanishing of Taylor series) on $\beta$.
  This is the {\bf dessin d'enfant}.
\end{definition}
Restricted by Riemann-Hurwitz, $\beta^{-1}(\infty)$ is not an independent degree of freedom, but is rather taken 1-to-1 to the faces in the bipartite graph, with the number of sides of the polygonal face being twice the ramification index.
One can succinctly record the valencies as
\[
\left[
  r_0(1), r_0(2), \ldots, r_0(B) \ \big| \
  r_1(1), r_1(2), \ldots, r_1(W) \ \big| \
  r_\infty(1), r_\infty(2), \ldots, r_\infty(I)
  \right]
\]
where $r_0(j)$ is the valency (ramification index) of the $j$-th black node, likewise $r_1(j)$, that for the $j$-th white node, and $r_\infty(j)$, half the number of sides to the $j$-th face.
This nomenclature is called the {\em passport} of the dessin; it does not uniquely determine it since one further needs the connectivity between the white/black nodes, but it is nevertheless an important quantity.
The {\em degree} of the Belyi map is the row-sum $d = \sum\limits_j r_0(j) = \sum\limits_j r_1(j) = \sum\limits_j r_\infty(j)$ and is the degree of $\beta$ as a rational function.
Finally, Riemann-Hurwitz demands that
\begin{equation}
2g-2 = d - (B+W+I)
\end{equation}
where $g$ is the genus of $\Sigma$ and $B,W,I$ are respectively the number of pre-images of $0,1,\infty$.
Of course, our focus will be on $g=1$, balanced (i.e., $B=W$) dessins, for which the total number of pre-images of $0,1,\infty$ is equal to the degree.
We emphasize that the dessin is {\em completely rigid}, fixing the coefficients in the defining polynomial of $\Sigma$ to be specific - and oftentimes horrendous - algebraic numbers with no room for complex parametres.
That is, given a bipartite graph on a Riemann surface, there is, up to coordinate change (birational transformation), a unique algebraic model and associated Belyi map which gives the graph as a dessin; even the slight change in the numerical coefficients in either will result in an utterly different graph.

In passing we mention that there is a purely combinatorial way to encode dessins - though, of course, in bypassing the analytics of Belyi maps, the algebro-number-theoretic information is somewhat lost.
Nevertheless, the ramification structure is captured completely, unlike the rather coarse encoding by passports.
One simply labels all edges, say $1,2, \ldots, d$, of the graph on $\Sigma$, and work within the symmetric group $\cS_n$ using the standard cycle notation for the elements.
Choose a {\em single} orientation, say clockwise, and for each black node, write the cycle $(e_1 e_2 \ldots e_i)$ going clockwise where $e_j$ is the integer label for the said edge.
Then form an element $\sigma_B \in \cS_d$ which is the product over all cycles for all black nides.
Similarly, do this for the white nodes (also going clockwise for the edges), forming $\sigma_W \in \cS_d$.
Finally, define $\sigma_\infty \in \cS_d$ such that
\begin{equation}
\sigma_B \sigma_W \sigma_\infty = Id \ ;
\end{equation}
this so-called {\bf permutation triple} captures the structure (valency and connectivity) of the dessin. Indeed, when $\sigma_\infty$ is itself expressed in terms of product of cycles - as every element of a symmetric group ultimately must - each cycle is associated to a face of the dessin, with its length being half the number of sides of the polygon.

Let us consider a concrete example \cite{Jejjala:2010vb}.
Consider the Belyi pair of elliptic curve $\Sigma$ and Belyi map $\beta$
\begin{equation}\label{BelyiC3}
y^2 = x^3 + 1 \ , \quad \beta = \frac12 ( 1 + y) \ .
\end{equation}
We see that the pre-image of 0 has $y = -1$, whence $(x,y) = \beta^{-1}(0) = (0,-1)$ on the elliptic curve $y^2 = x^3 + 1$. Choosing local coordinates $(x,y) = (0 + \epsilon, -1 + \delta)$ for infinitesimals $(\epsilon, \delta)$, we have that $-2 \delta = \epsilon^3$. Therefore locally the map is $\beta = \frac{\delta}{2} \sim -\frac{\epsilon^3}{4}$ and thus the ramification index for this single pre-image of 0 is 3.
Similarly, the pre-image of 1 has $y = +1$, i.e., $(x,y) = (0,1)$ is the single pre-image of 1, where local coordinates can be chosen as $(0 + \epsilon, 1 + \delta)$, so that $\beta \sim \frac{\epsilon^3}{4}$. Hence, the ramification index is also 3 for this single pre-image of 1.
Finally, $(\infty, \infty)$ is the pre-image of $\infty$ where the local coordinates $(\epsilon^{-1}, \delta^{-1})$ can be chosen so that $(x,y) \sim (\epsilon^{-2}, \epsilon^{-3})$. Hence, the ramification index at $\infty$ is also 3.
Thus, the passport is $[3 \big| 3 \big| 3]$ and the maps is degree 3.
In summary
\[
\begin{array}{|c|c|c|c|c|c|}
\hline
\mathbb{T}^2: y^2 = x^3 + 1 & \stackrel{\beta=\frac12(1+y)}{\longrightarrow} & \mathbb{P}^1 & \mbox{Local Coordinates on }\mathbb{T}^2 & \mbox{Ram.~Index}(\beta) \\ \hline\hline
(0,-1) & \sim -\frac14{\epsilon^3} & 0 & (x,y) \sim (\epsilon,-1-\frac12\epsilon^3) & 3 \\ \hline
(0,1) & \sim \frac14{\epsilon^3} & 1 & (x,y) \sim (\epsilon,1+\frac12\epsilon^3) & 3 \\ \hline
(\infty,\infty) & \sim \frac12 \epsilon^{-3} & \infty & (x,y) \sim (\epsilon^{-2},\epsilon^{-3}) & 3 \\ \hline
\end{array}
\]
Therefore, what we have is a bipartite graph on $T^2$ specified by the elliptic curve in \eqref{BelyiC3}, with a single pair of black/white nodes, each trivalent.
This is none other than \eqref{c3dimer}, but with completely rigid shape:
the j-invariant for this elliptic curve is $j = 0$, exactly that of the a-maximized isoradial embedding in \eqref{aMaxC3}.

It must be stressed that whilst drawing the dessin knowing the Belyi pair is a fairly straight-forward matter, knowing a bipartite graph on a Riemann surface and reconstructing the Belyi pair is a tremendously difficult task and no general method exists.
There is an entire field devoted to finding explicit (q.v.~\cite{catalog,hv,vk,Vidunas:2016xun}) Belyi maps and hyper-elliptic models.
Part of the power of Grothendieck's original insight was that one did not need to know the explicit Belyi maps in order to study dessins and that one of the central objects of number theory, the absolute Galois group $Gal(\overline{\IQ}/\IQ)$, acts faithfully on the space of dessins.
It is perhaps ironic that now, for the sake of both physics and mathematics, one is in need of these explicit maps.

Continuing with our examples, viz., $\IF_0(I)$ and $\IF_0(II)$, we can likewise check \cite{Hanany:2011ra} that
\begin{equation}
  \begin{array}{|c|c|c|} \hline
    & F_0(I) & F_0(II) \\ \hline
    \text{Elliptic Curve} & y^2=x^3-x & y^2=x^3-x \\ \hline
    \text{Belyi Map} & i\frac{(i+x)^4}{8\,x\,(1-x^2)} &
    i\,\frac{(x^2-(-1)^{1/3})^3}{3\,\sqrt{3}\,x^2\,(x^2-1)} \\ \hline
    \text{Passport} &
    [4,4 \ \big| \ 4,4 \ \big| \ 2,2,2,2] 
    &
    [3,3,3,3 \ \big| \ 3,3,3,3 \ \big| \ 2,2,4,4]  \\ \hline
  \end{array}
\end{equation}
We note that the form of the elliptic curves are the same.
Moreover, in this example, the j-invaraint of both $\IF_0$ theories are equal to 1, matching the results from \eqref{F0IR} and \eqref{F0IR}.

One might be tempted to conclude that $j(\tau)$ of the dessin is equal to the $j(\tau)$ of the a-maximized isoradial bipartite graph and that, parallel to Theorem \ref{clusterj}, $j(\tau)$ is also an invariant for dessins under urban renewal.
Sadly, both such conjectures \cite{Jejjala:2010vb} were later shown to have counter-examples \cite{Hanany:2011ra,He:2012xw,Vidunas:2016xun}.
It remains, however, an interesting question as to how one may modify, for instance, the isoradial embedding or the a-function, so that this identification of the elliptic curve of the dessin matches that from the quiver side.

Conversely, the dessin community usually thinks in depth about the orbits under which classes of dessins fall via Galois actions on the algebraic numbers involved.
Orbits under urban renewal, whose incipience came from cluster algebras and quantum field theory, would be a rather novel object to ponder.
Now, it is certainly not true that urban renewal stays in the same Galois orbit \cite{Vidunas:2016xun}.
However, inspired by part 2 of Theorem \eqref{clusterj} on the invariance of transcendence degree in the a-maximization case, we should look into the matter further.

\section{Conclusions and Prospectus}\label{s:conc}
We have embarked upon an extraordinary journey, from the representation theory of quivers, to the geometry of Calabi-Yau varieties, from the bipartite graphs of  dimer models to the intricacies of the Grothendieck's dessin d'enfants.
This is a journey from algebra to geometry to number theory, all under the auspices of the theoretical physics of gauge/string theories.
Whilst we hope to have tantalized the reader, it is self-evident that we have only touched upon the surface, and that there remain a multitude of open problems awaiting further investigation.

As for the VMS, the question arises as to whether Theorem \ref{seiberg} and its Corollary \ref{equalM} exhaust the possible isomorphisms as affine varieties.
One might expect that there should be dualities generalizing Seiberg duality/cluster mutation/urban renewal, but what are they and what are the implications for quivers and for the quantum field theory thereby encoded?
As we saw that the quiver variety $\cM(\cQ)$ being Calabi-Yau is only an extremely special case, how does our story persist for arbitrary affine varieties? How does the plethystic programme on the Hilbert series, which in the Calabi-Yau case detects the Gorenstein property through palindromy, aid such general treatment?
Indeed, what are the geometrical consequences of Observations \ref{plog} in disentangling generators and syzygies, in a similar spirit to how Molien series intertwine generators and relations for finite groups.

Now, the bipartite structure of $\cB$ also seems inextricably linked to the toric nature of the VMS, and the underlying torus of the brane-tiling, to superconformality.
The planar nature of the tiling is also a reflection of the Calabi-Yau nature of $\cM(\cQ)$, which ensures that the endpoints of lattice vectors generating the toric cone be co-planar.
The ubiquity of bipartite structure, from cluster algebras to the positive Grassmannian to the amplitudes of super-Yang-Mills theories, has permeated mathematics and physics of late.
It is clearly important to understand it in more depth and to also study the results of proceeding beyond (balanced) bipartite graphs.

The a-maximization and its geometric interpretation of isoradial embedding of the tiling subject to the extremization of a cubic in accordance with Theorem \ref{amax} places a strigent constraint on the shape of the bipartite graph.
Interestingly, this shape is invariant under urban renewal (cf.~Theorem \ref{clusterj}) and the transcendence degree of the R-charges over $\IQ$ is also unchanged.
As we venture into number theory, the manifestation of $\cB$ as dessin d'enfant is almost unavoidable.
The rigidification of the graph to specific algebraic numbers is a remarkable fact.
However, unlike the a-maximization case where the algebraic numbers to which the R-charges freeze have a clear meaning in the dual Calabi-Yau geometry as normalized Sasaki-Einstein volume, here they are mysterious.

Indeed, what is urban renewal, central to the quiver $\cQ$,  in the theory of dessins?
Conversely, what is Galois conjugation, key to the dessin $\cB$, in the quiver?
These, and many more questions, beckon our continued exploration.

\newpage 
\section*{Acknowledgments}
{\it Catharinae Sanctae Alexandriae adque Majorem Dei Gloriam}\\
As I express my gratitude to Professors Ji, Papadopoulos, Schneps, \& Su for organizing the excellent workshop on ``Grothendieck-Teichm\"uller theories'' at the Chern Institute, Nankai, which engendered much lively and helpful conversations, and to which they kindly invited me to contribute, I would also like to lend this opportunity to thank my collaborators and friends, with whom I have had countless enjoyable discussions on the material presented here: foremost, Professor Amihay Hanany, paternal in my student days, and now avuncular, but always perspicacious in thought and prophetic in words.

I toast to my brethren Sergio Benvenuti, Bo Feng, Davide Forcella, Sebastian Franco, Jon Hauenstein, Vishnu Jejjala, Kris Kennaway, Alastair King, Dhagash Mehta, Noppadol Mekareeya, Jurgis Pasukonis, Sanjaye Ramgoolam, Diego Rodriguez-Gomez, Rak-Kyeong Seong, James Sparks, Angel Uranga, Raimundas Vidunas, and Alberto Zaffaroni for their cheerful company, and to my talented students and postdocs Nessi Benishti, Sownak Bose, James Gundry, Mark van Loon, Malte Probst, James Read, Chuang Sun, Yan Xiao, Zhi Hu, and Da Zhou, for their constantly re-instilling in me a sense of youth.

I am grateful to Professor S.-T.~Yau and Professor Cumrum Vafa for their perpetual guidance, and especially to Professor John McKay for his endless wisdom and sagacious warmth, to whose penetrating insights I pay awe-struck homage.

I am indebted to the Science and Technology Facilities Council, UK, for grant ST/J00037X/1, the Chinese Ministry of Education, for a ChangJiang Chair Professorship at NanKai University, and the city of Tian-Jin for a Qian-Ren Award.
Above all, I thank Merton College, Oxford for continuing to provide a quiet corner of Paradise for musings and contemplations.


\end{document}